\documentclass[11pt,twoside]{article}

\usepackage{amssymb,amsmath,amsfonts,amsthm,color,mathrsfs}
\usepackage[Symbol]{upgreek}
\usepackage{txfonts}
\usepackage[nottoc,notlot,notlof]{tocbibind}
\usepackage[active]{srcltx}
\usepackage{hyperref}
\usepackage{citeref}
\usepackage{mathrsfs}
\usepackage[top=1in, bottom=1in, left=1.2in, right=1.2in]{geometry}

\theoremstyle{plain}

\newtheorem{thm}{Theorem}[section]
\newtheorem{lem}[thm]{Lemma}
\newtheorem{prop}[thm]{Proposition}

\newtheorem{defn}[thm]{Definition}

\newtheorem{con}[thm]{Conclusion}
\numberwithin{equation}{section}

\allowdisplaybreaks
\begin{document}

\arraycolsep=1pt
\author{ Zaihui Gan$^{\sharp}$  ~~~~ Yue Wang  \\
\\
{ Center for Applied Mathematics, Tianjin University,  Tianjin 300072, China}}

\title{ { Existence and Instability of Standing Wave for the Two-wave Model with Quadratic Interaction }
 \footnotetext{\hspace{-0.35cm}
{$\sharp$: Corresponding authors: Zaihui Gan}}
 \footnotetext{\hspace{-0.35cm}
{Emails: ganzaihui2008cn@tju.edu.cn (Zaihui Gan), wangy2017@tju.edu.cn (Yue Wang).
 }}
 }
\date{}
\maketitle

\begin{center}
\begin{minipage}{13cm}
{\noindent{\bf Abstract.} In this paper, we establish the existence and instability of standing wave for a system of nonlinear Schr\"{o}dinger equations arising in the two-wave model with quadratic interaction in higher space dimensions under mass resonance conditions. Here, we eliminate the limitation for the relationship between complex constants $a_{1}$ and $a_{2}$ given in \cite{HOT}, and consider arbitrary real positive constants $a_{1}$ and $a_{2}$.  First of all, according to the conservation identities for mass and energy, using the so-called virial type estimate, we obtain that the solution for the Cauchy problem under consideration blows up in finite time in $H^{1}(\mathbb{R}^{N})\times H^{1}(\mathbb{R}^{N})$ with space dimension $N\geq 4$. Next, for space dimension $N$ with $4<N<6$, we establish the existence of the ground state solution for the elliptic equations corresponding to the nonlinear Schr\"{o}dinger equations under the frequency and mass resonance by adopting variational method, and further achieve the exponential decay at infinity for the ground state. This implies the existence of standing wave for the nonlinear Schr\"{o}dinger equaitons under consideration. Finally, by defining another constrained minimizing problems for a pair of complex-valued functions, a suitable manifold, referring to the characterization of the standing wave, making appropriate scaling and adopting virial type estimate, we attain the instability of the standing wave for the equations under frequency and mass resonance in space dimension $N$ with $4<N<6$ by virtue of the conservations of mass and energy. Here, we adopt the equivalence of two constrained minimizing problems defined for pairs of complex-valued and real-valued functions $(u,v)$, respectively, when $(u,v)$ is a pair of real-valued functions.
}
\end{minipage}
\end{center}
\vspace{0.2cm}
{\bf Statements and Declarations:} No conflict of interest exists in the submission of this manuscript. No data was used for the research described in this manuscript.

%
%
%
%
%

\section{Introduction}
\indent{\quad}We consider in this paper the nonlinear Schr\"{o}dinger equations in the two-wave model with quadratic interaction

\begin{equation}\label{1.1}
\left\{~~
\begin{array}{ll}
&i\phi_t+\dfrac{1}{2m_1}\Delta\phi=a_1\psi\overline{\phi},\quad t>0,\quad x\in\mathbb{R}^N,\quad(1.1-a)\\[0.5cm]
&i\psi_t+\dfrac{1}{2m_2}\Delta\psi=a_2\phi^2,~\quad t>0,\quad x\in\mathbb{R}^N,\quad(1.1-b)
\end{array}
\right.
\end{equation}
where $4\leq N< 6$, $\phi=\phi(t,x)$ and $\psi=\psi(t,x)$ are complex-valued functions of $(t,x)\in \mathbb{R}^+\times\mathbb{R}^N$; $\Delta$ is the Laplacian in $\mathbb{R}^N$; $m_1,m_2,a_1$ and $a_2$ are positive constants, and $\overline{\phi}$ is the complex conjugate of $\phi$.\\
\indent Note that the modulated wave functions $(\phi_c,\psi_c)=(e^{itm_1c^2}\phi,e^{itm_2c^2}\psi)$ satisfy
\begin{equation}\label{1.2}
\left\{
\begin{array}{ll}
&\dfrac{1}{2c^2m_1}\partial_t^2\phi_c-i\partial_t\phi_c-\dfrac{1}{2m_1}\Delta\phi_c
=-e^{itc^2(2m_1-m_2)}a_1\psi_c\overline{\phi}_c, \\[0.5cm]
&\dfrac{1}{2c^2m_2}\partial_t^2\psi_c-i\partial_t\psi_c-\dfrac{1}{2m_2}\Delta\psi_c
=-e^{itc^2(m_2-2m_1)}a_2\phi_c^2,
\end{array}
\right.
\end{equation}
which formally leads to \eqref{1.1} as the speed of light $c$ tends to infinity under the  mass resonance condition:
\begin{equation}\label{1.3}
m_2=2m_1.
\end{equation}
On the other hand, by \eqref{1.2} one knows that $(\phi,\psi)$ verifies the nonlinear Klein-Gordon equations
\begin{equation}\label{1.4}
\left\{
\begin{array}{ll}
&\dfrac{1}{2c^2m_1}\partial_t^2\phi-\dfrac{1}{2m_1}\Delta\phi+\dfrac{m_1c^2}{2}\phi=-a_1\psi\overline{\phi},\\\\
&\dfrac{1}{2c^2m_2}\partial_t^2\psi-\dfrac{1}{2m_2}\Delta\psi+\dfrac{m_2c^2}{2}\psi=-a_2\phi^2.
\end{array}
\right.
\end{equation}
Thus the system \eqref{1.1} can be regards as the non-relativistic limit of the nonlinear Klein-Gordon equations \eqref{1.4} (see also \cite{HOT}).\\
\indent A number of authors (\cite{KOS,MNO1,MNO2,Suna1,Tsut1}) have established some results concerning the non-relativistic limit for the nonlinear Klein-Gordon equations. In view of the relationship between the nonlinear Klein-Gordon equations (1.4) and the nonlinear Schrodinger equations (1.1), in \cite{HOT}, the authors established some results for the system \eqref{1.1} including local existence and global existence of $L^2-$ solution and $H^1-$ solution, Galilei invariance of local and global solutions under mass resonance $(m_2=2m_1)$, non-existence of global solutions with negative energy under mass resonance $(m_2=2m_1)$ and the existence of ground states. These results are valid under the conditions $a_1=c\overline{a}_2$ for a constant $c\in\mathbb{R}\backslash\{0\}$, and $a_1$, $a_2$ are complex constants.\\
\indent In the present paper for the equations \eqref{1.1}, however, we eliminate the assumption that $a_1=c\overline{a}_2$ for some $c\in\mathbb{R}\backslash\{0\}$ given in \cite{HOT} and consider the cases $a_1$ and $a_2$ are real constants with arbitrary $a_1>0$ and $a_2>0$. In particular, we consider first here the stationary equations corresponding to equations \eqref{1.1} and then establish the instability of the standing wave for \eqref{1.1}. \\
\indent Let $\omega_1>0$ and $\omega_2=2\omega_1>0$ be two real constants, and $(\phi(t,x),\psi(t,x))=(e^{i\omega_1t}u(x),e^{i\omega_2t}v(x))$ for a pair of real functions $(u,v)=(u(x),v(x))$. Then $(u(x),v(x))$ is the ground state solution to the following equations under the case $\omega_2=2\omega_1$:
\begin{equation}\label{1.5}
\left\{
\begin{array}{ll}
&-\dfrac{1}{2m_1}\Delta u(x)+\omega_1u(x)=-a_1v(x)u(x),\quad (1.5a)\\\\
&-\dfrac{1}{2m_2}\Delta v(x)+\omega_2v(x)=-a_2u^2(x),~\qquad (1.5b)
\end{array}
\right.
\end{equation}
and $(u(x),v(x))\in H^1(\mathbb{R}^N)\times H^1(\mathbb{R}^N)\backslash\{(0,0)\}$. Thus, $(\phi(t,x),\psi(t,x))=(e^{i\omega_1t}u(x),e^{i\omega_2t}v(x))$ is a standing wave solution of the evolution equations \eqref{1.1}. On the other hand,from the physical viewpoint, the ground state solution of \eqref{1.5} plays a key role. Recalling that the definition of the ground state for a single elliptic equation mentioned in Cazenave \cite{Caze1}, we can define the ground state for the elliptic system \eqref{1.5}. A nontrivial solution $(u,v)$ of \eqref{1.5} is named a ground state if it is of minimum action among all solutions of \eqref{1.5}. Specifically,  $(u,v)$ satisfies
 $$S(u,v)\leq S(\omega,\xi)$$
for any solution $(\omega,\xi)$ of \eqref{1.5}, where the action $\mathcal{S}(u,v)$ for a pair of real functions $(u,v)$ is defined by
\begin{equation}\label{1.6}
\left.
\begin{array}{ll}
S(u,v)=&\displaystyle\dfrac{a_2}{2m_1}\int_{\mathbb{R}^N}|\nabla u|^2dx+\dfrac{a_1}{4m_2}\int_{\mathbb{R}^N}|\nabla v|^2dx\\\\
&\displaystyle+a_2\omega_1\int_{\mathbb{R}^N}|u|^2dx+\frac{a_1}{2}\omega_2\int_{\mathbb{R}^N}|v|^2dx+
a_1a_2\int_{\mathbb{R}^N}vu^2dx.
\end{array}
\right.
\end{equation}
In the present paper, we consider three aspects for the nonlinear Schr\"{o}dinger equations \eqref{1.1}: one is blowup in finite time for space dimension $ 4\leq N<6$, the others are the existence and instability of the standing wave for space dimension $4< N<6 $ . For obtaining that the solution to the Cauchy problem for \eqref{1.1} blows up in finite time in $H^{1}(\mathbb{R}^N)\times H^{1}(\mathbb{R}^N)$, we need to establish the so-called virial type estimates and to utilize the modified conservation identities for mass and energy. On the other hand, the key to attain the existence of standing wave for  equations \eqref{1.1} with space dimension $4< N<6$ is to attain the existence of the ground state solution for equations \eqref{1.5}, where the ground state solution $(\xi,-\eta)$ is dependent of $|x|$ alone and has an exponential decay at infinity. Furthermore, for exploring the instability of the standing wave for \eqref{1.1} with space dimension $ 4< N<6$, we must be concerned with the characterization of the standing waves for \eqref{1.1} with minimal action $S(\xi,-\eta)$, and refer to the exponential decay at infinity of the ground state $(\xi,-\eta)$ of \eqref{1.5}. Here, $(\xi,-\eta)$ will be decided in Section 4.\\
 \indent Throughout this paper, we will adopt a type of variational method initially proposed in \cite{Bere-Lions1,Bere-Lions2}; the main thing for the variational method is to define suitable functionals, manifolds and constrained minimizing problem.
Here we define the functional $Q(u,v)$ for a pair of real-valued functions $(u,v)$ by
\begin{equation}\label{1.7}
\left.
\begin{array}{ll}
Q(u,v)=&\displaystyle\dfrac{a_2}{m_1}\int_{\mathbb{R}^N}|\nabla u|^2dx+\dfrac{a_1}{2m_2}\int_{\mathbb{R}^N}|\nabla v|^2dx
+\dfrac{N}{2}a_1a_2\int_{\mathbb{R}^N}v u ^2dx,
\end{array}
\right.
\end{equation}
the manifold $\mathcal{M}$ as
\begin{equation}\label{1.8}
\left.
\begin{array}{ll}
\mathcal{M}:=\left\{(u,v)\in H^1(\mathbb{R}^N)\times H^1(\mathbb{R}^N)\backslash\{(0,0)\},Q(u,v)=0\right\},
\end{array}
\right.
\end{equation}
and the constrained minimizing problem as
\begin{equation}\label{1.9}
\left.
\begin{array}{ll}
K=\inf\limits_{(u,v)\in\mathcal{M}}S(u,v).
\end{array}
\right.
\end{equation}
Throughout this paper, we make the following assumptions on these real coefficients $m_{1},~m_{2},~a_{1},~a_{2},$ in equations \eqref{1.1}:
\begin{equation}\label{1.10}
\left.
\begin{array}{ll}
m_1>0,m_2>0,a_1>0,a_2>0.
\end{array}
\right.
\end{equation}
\indent This paper is organized as follows. Section 2 is devoted to establishing some basic estimates and to collecting basic lemmas which will be used in the subsequent sections. In section 3, we establish finite time blow up for solution of \eqref{1.1} in terms of the initial data $(\phi_0,\psi_0)$ and the initial energy. In section 4, we show the existence of the ground state solution and its exponential decay at infinity. In section 5, we justify the instability of the standing waves in view of these conclusions given in the section 2, section 3 and section 4.
\section{Preliminaries}
We impose the initial data on \eqref{1.1} as follows:
\begin{equation}\label{2.1}
\left.
\begin{array}{ll}
\phi(0,x)=\phi_0(x),\quad\psi(0,x)=\psi_0(x),\quad x\in\mathbb{R}^N,
\end{array}
\right.
\end{equation}
 where $(\phi_0(x),\psi_0(x))$ are given functions in $H^1(\mathbb{R}^N)\times H^1(\mathbb{R}^N)$.
 From the results of Ginibre and Velo \cite{GV1}, Cazenave\cite{Caze1} as well as Hayashi, Ozawa and Tanakai \cite{HOT}, \eqref{1.1}-\eqref{2.1} is locally well-posed in $H^1(\mathbb{R}^N)\times H^1(\mathbb{R}^N)$, and thus for any $(\phi_0,\psi_0)\in H^1(\mathbb{R}^N)\times H^1(\mathbb{R}^N)$, there exists a unique solution $(\phi,\psi)$ to the Cauchy problem \eqref{1.1}-\eqref{2.1} in $C\left([0,T);H^1(\mathbb{R}^N)\times H^1(\mathbb{R}^N)\right)$ defined on a maximal time interval $T_{max}(\phi_0,\psi_0)$, either $T=+\infty$ or $T<+\infty$ and
 $$\lim\limits_{t\rightarrow T_{max}^-(\phi_0,\psi_0)}\left(
 \|\phi\|_{H^1(\mathbb{R}^N)}+\|\psi\|_{H^1(\mathbb{R}^N)}\right)=+\infty.$$
In addition, the Cauchy problem \eqref{1.1}-\eqref{2.1} admits the following conservation laws in the energy space $H^1(\mathbb{R}^N)\times H^1(\mathbb{R}^N)$.
\begin{lem}\label{Lemma 2.1}
Suppose that $(\phi,\psi)\in C\big([0,T);H^1(\mathbb{R}^N)\times H^1(\mathbb{R}^N)\big)$ is a solution of the Cauchy problem \eqref{1.1}-\eqref{2.1}. Then the total mass, total energy and the total momentum are conserved for all $t\geq0$:
 \\[0.3cm]
{\bf $L^2$-norm (Mass):}
\begin{equation}\label{2.2}
\left.
\begin{array}{ll} \displaystyle\int_{\mathbb{R}^N}\left(a_2|\phi(t,x)|^2+a_1|\psi(t,x)|^2\right)dx
=\int_{\mathbb{R}^N}\left(a_2|\phi_0(x)|^2+a_1|\psi_0(x)|^2\right)dx;
\end{array}
\right.
\end{equation}
{\bf Energy:}
\begin{equation}\label{2.3}
\left.
\begin{array}{ll}
 E(\phi(t,x),\psi(t,x))=E(\phi_0(x),\psi_0(x)),
\end{array}
\right.
\end{equation}
where $E\left(\phi(t,x),\psi(t,x)\right)$ is defined as
\begin{equation*}\label{2.3a}
\left.
\begin{array}{ll}
  E(\phi(t,x),\psi(t,x)) &\displaystyle=\dfrac{a_2}{2m_1}\int_{\mathbb{R}^N}|\nabla \phi(t,x)|^2dx+\dfrac{a_1}{4m_2}\int_{\mathbb{R}^N}|\nabla \psi(t,x)|^2dx
  \\[0.5cm]
 &\displaystyle\quad
+a_1a_2Re\int_{\mathbb{R}^N}\psi(t,x)\overline{\phi}(t,x)^2dx;
\end{array}
\right.\eqno(2.3a)
\end{equation*}
{\bf Momentum:}
\begin{equation}\label{2.4}
\left.
\begin{array}{ll}
 &\displaystyle a_2Im\left(\int_{\mathbb{R}^N}\nabla \phi(t,x)\overline{\phi}(t,x)dx\Big)+
 \dfrac{1}{2}a_1Im\Big(\int_{\mathbb{R}^N}\nabla \psi(t,x)\overline{\psi}(t,x)dx\right)\\\\
 &\quad= \displaystyle a_2Im\left(\int_{\mathbb{R}^N}\nabla \phi_0(x)\overline{\phi}_0(x)dx\right)+
 \dfrac{1}{2}a_1Im\left(\int_{\mathbb{R}^N}\nabla \psi_0(x)\overline{\psi}_0(x)dx\right).
\end{array}
\right.
\end{equation}
\end{lem}
\begin{proof}
Multiplying $(1.1a)$ by $a_2\overline{\phi}$ and $(1.1b)$ by $a_1\overline{\psi}$, integrating over $\mathbb{R}^N$ and taking the imaginary part for the resulting equations, we obtain formally
\begin{equation*}
\left.
\begin{array}{ll}
&\displaystyle\dfrac{d}{dt}\left(\dfrac{1}{2}
\int_{\mathbb{R}^N}\left(a_2|\phi(t,x)|^2+a_1|\psi(t,x)|^2\right)dx\right)
\\[0.5cm]
 & \quad= \displaystyle Im\left(a_1a_2\int_{\mathbb{R}^N}\psi(t,x)\overline{\phi}^2(t,x)dx\right)
 +Im\left(a_1a_2\int_{\mathbb{R}^N}\overline{\psi}(t,x)\phi^2(t,x)dx\right)
 \\[0.5cm]
 &\quad=0,
\end{array}
\right.
\end{equation*}
which implies the conservation of mass
$$\int_{\mathbb{R}^N}\left(a_2|\phi(t,x)|^2+a_1|\psi(t,x)|^2\right)dx
=\int_{\mathbb{R}^N}\left(a_2|\phi_0(x)|^2+a_1|\psi_0(x)|^2\right)dx
.$$
Next, multiplying $(1.1a)$ by $2a_2\overline{\phi}_t$ and $(1.1b)$ by $a_1\overline{\psi}_t$, integrating over $\mathbb{R}^N$ and taking the real part, we obtain
\begin{equation*}
\left.
\begin{array}{ll}
  \displaystyle\dfrac{d}{dt}\left[\dfrac{a_2}{2m_1}\int_{\mathbb{R}^N}|\nabla \phi(t,x)|^2dx+\dfrac{a_1}{4m_2}\int_{\mathbb{R}^N}|\nabla \psi(t,x)|^2dx
+a_1a_2Re\int_{\mathbb{R}^N}\psi(t,x)\overline{\phi}(t,x)^2dx\right]=0,
\end{array}
\right.
\end{equation*}
which yields formally the conservation of energy
$$E\left(\phi(t,x),\psi(t,x)\right)=E\left(\phi_0(x),\psi_0(x)\right),$$
where the energy $E$ is defined by $(2.3a)$.
Finally, multiplying $(1.1a)$ by $2a_2\nabla\overline{\phi}$ and $(1.1b)$ by $a_1\nabla\overline{\psi}$, integrating over $\mathbb{R}^N$ and taking the real part, we obtain
\begin{equation*}
\left.
\begin{array}{ll}
 &\displaystyle Re\int_{\mathbb{R}^N}2 a_2i\phi_t \nabla\overline{\phi}dx+Re\int_{\mathbb{R}^N}a_1i\psi_t\nabla\overline{\psi}dx
 \\[0.5cm]
 &\displaystyle\quad=-Re\int_{\mathbb{R}^N} \dfrac{ a_2}{ m_1}\nabla \overline{\phi}\Delta\phi dx-Re\int_{\mathbb{R}^N} \dfrac{a_1}{2m_2}\nabla \overline{\psi}\Delta\psi dx
 \\[0.5cm]
 &\displaystyle\quad\quad+Re\int_{\mathbb{R}^N}2a_1a_2\psi
 \overline{\phi}\nabla\overline{\phi}dx+Re\int_{\mathbb{R}^N}a_1a_2\phi^2\nabla\overline{\psi}dx.
 \end{array}
\right.
\end{equation*}
Direct calculation gives
$$\dfrac{d}{dt}\left[ a_2Im\left(\int_{\mathbb{R}^N}\nabla \phi(t,x)\overline{\phi}(t,x)dx\right)+
 \frac{1}{2}a_1Im\left(\int_{\mathbb{R}^N}\nabla \psi(t,x)\overline{\psi}(t,x)dx\right)\right]=0,$$
 which yields formally the conservation of momentum
\begin{equation*}
\left.
\begin{array}{ll}
 &\displaystyle a_2Im\left(\int_{\mathbb{R}^N}\nabla \phi\overline{\phi}dx\right)+
 \dfrac{1}{2}a_1Im\left(\int_{\mathbb{R}^N}\nabla \psi\overline{\psi}dx\right)\\\\
  &\displaystyle\quad=a_2Im\left(\int_{\mathbb{R}^N}\nabla \phi_0\overline{\phi}_0dx\right)+
 \dfrac{1}{2}a_1Im\left(\int_{\mathbb{R}^N}\nabla \psi_0\overline{\psi}_0dx\right).
\end{array}
\right.
\end{equation*}
\end{proof}
Let
\begin{equation}\label{2.5}
\displaystyle \Sigma:= \left\{u:~~u\in H^1(\mathbb{R}^N),~~ xu\in L^2(\mathbb{R}^N)\right\}.
\end{equation}
We then establish a kind of virial type estimate for the Cauchy problem \eqref{1.1}-\eqref{2.1}, which is helpful for exploring the instability of standing waves.
\begin{lem}\label{Lemma 2.2}
Assume that $(\phi_0, \psi_0)\in \Sigma\times\Sigma$ . Let $(\phi,\psi)\in C\big([0,T);H^1(\mathbb{R}^N)\times H^1(\mathbb{R}^N)\big)$ be the corresponding solution of the cauchy problem \eqref{1.1}-\eqref{2.1} on $[0,T)$.\\
\indent Put
\begin{equation}\label{2.6}
\left.
\begin{array}{ll}
\displaystyle G(t)=\int_{\mathbb{R}^N}|x|^2\left(a_2|\phi|^2+a_1|\psi|^2\right)dx,
\end{array}
\right.
\end{equation}
then it follows that $$\left(|x|\phi(t,.),|x|\psi(t,.)\right)\in C\left((-T_{min},T_{max}),L^2(\mathbb{R}^N) \right)\times C\left((-T_{min},T_{max}),L^2(\mathbb{R}^N)\right).$$
Moreover, there hold
\begin{equation}\label{2.7}
\left.
\begin{array}{ll}
\displaystyle\dfrac{d}{dt}G(t)=4Im\int_{\mathbb{R}^N}
\left(\dfrac{a_2}{2m_1}x\overline{\phi}\nabla\phi
+\dfrac{a_1}{2m_2}x\overline{\psi}\nabla\psi\right)dx,
\end{array}
\right.
\end{equation}
and
\begin{equation}\label{2.8}
\left.
\begin{array}{ll}
\displaystyle\dfrac{d^2}{dt^2}G(t)&\displaystyle=\dfrac{2a_2}{m_1^2}
\int_{\mathbb{R}^N}|\nabla\phi|^2dx+\dfrac{2a_1}{m_2^2}
\int_{\mathbb{R}^N}|\nabla\psi|^2dx
\\[0.5cm]
&\quad\displaystyle+\dfrac{2a_1a_2N}{m_2}Re\int_{\mathbb{R}^N}
\overline{\psi}\phi^2dx+2a_1a_2\left(\dfrac{2}{m_2}-\dfrac{1}{m_1}\right)
Re\int_{\mathbb{R}^N}x\phi^2\nabla\overline{\psi}dx
\\[0.5cm]
&\displaystyle=\dfrac{2a_2}{m_1^2}\int_{\mathbb{R}^N}|\nabla\phi|^2dx
+\dfrac{2a_1}{m_2^2}\int_{\mathbb{R}^N}|\nabla\psi|^2dx
\\[0.5cm]
&\quad\displaystyle+2a_1a_2N\left(\dfrac{1}{m_1}-\dfrac{1}{m_2}\right)
Re\int_{\mathbb{R}^N}\overline{\psi}\phi^2dx+2a_1a_2\left(\frac{1}{m_1}
-\frac{2}{m_2}\right)Re\int_{\mathbb{R}^N}x\nabla\phi^2\overline{\psi}dx.
\end{array}
\right.
\end{equation}
\end{lem}
\begin{proof}
Since $(\phi,\psi)\in C\left([0,T);H^1(\mathbb{R}^N)\times H^1(\mathbb{R}^N)\right)$ is a solution of the Cauchy problem \eqref{1.1}-\eqref{2.1}, applying the results of Ginibre and Velo \cite{GV1}, from $(|x|\phi_0,|x|\psi_0)\in L^2(\mathbb{R}^N)\times L^2(\mathbb{R}^N)$, it follows that $(|x|\phi,|x|\psi)\in L^2(\mathbb{R}^N)\times L^2(\mathbb{R}^N)$. This implies that $G(t)$ given by \eqref{2.6} is well defined on $[0,T)$.\\
\indent Differentiating \eqref{2.6} with respect to $t$ yields
\begin{equation*}\label{(2.8-1)}
\left.
\begin{array}{ll}
 G'(t)&=\displaystyle\int_{\mathbb{R}^N}|x|^2
\left[a_2\left(\phi\overline{\phi}_t+\phi_t\overline{\phi}\right)
+a_1\left(\psi\overline{\psi}_t+\psi_t\overline{\psi}\right)\right]dx
\\[0.5cm]
&\displaystyle=2Re\int_{\mathbb{R}^N}|x|^2\left(a_2\phi\overline{\phi}_t+a_1\psi
\overline{\psi}_t\right)dx
\\[0.5cm]
&\displaystyle=-4Im\int_{\mathbb{R}^N}\left(\dfrac{a_2}{2m_1}x\phi\nabla
\overline{\phi}+\dfrac{a_1}{2m_2}x\psi\nabla\overline{\psi}\right)dx.
\end{array}
\right.\eqno(2.8-1)
\end{equation*}

Note that
\begin{equation*}
\left\{~~
\begin{array}{ll}
\displaystyle\phi_t=\dfrac{i}{2m_1}\Delta\phi-ia_1\psi\overline{\phi},\quad \overline{\phi}_t=-\dfrac{i}{2m_1}\Delta\overline{\phi}+ia_1\overline{\psi}\phi,
\\[0.5cm]
\displaystyle\psi_t=\dfrac{i}{2m_2}\Delta\psi-ia_2\phi^2,\quad \overline{\psi}_t=-\dfrac{i}{2m_2}\Delta\overline{\psi}+ia_2\overline{\phi}^2,
\\[0.5cm]
\displaystyle Re\int_{\mathbb{R}^N}x\Delta\overline{\phi}\nabla\phi dx=\dfrac{N-2}{2}\int_{\mathbb{R}^N}|\nabla\phi|^2dx,\quad Re\int_{\mathbb{R}^N}x\Delta\overline{\psi}\nabla\psi dx=\dfrac{N-2}{2}\int_{\mathbb{R}^N}|\nabla\psi|^2dx,
\\[0.5cm]
\displaystyle Re\int_{\mathbb{R}^N}x\phi\nabla\Delta\overline{\phi} dx=\dfrac{N+2}{2}\int_{\mathbb{R}^N}|\nabla\phi|^2dx,\quad Re\int_{\mathbb{R}^N}x\psi\nabla\Delta\overline{\psi} dx=\dfrac{N+2}{2}\int_{\mathbb{R}^N}|\nabla\psi|^2dx,
\\[0.5cm]
\displaystyle Re\int_{\mathbb{R}^N}\left(-x\psi\overline{\phi}\nabla\overline{\phi}\right)dx
+Re\int_{\mathbb{R}^N}x\left(\phi^2\nabla\overline{\psi}+\phi\overline{\psi}
\nabla\phi\right)dx
=Re\int_{\mathbb{R}^N}x\phi^2\nabla\overline{\psi}dx,
\end{array}
\right.\eqno(2.8-2)
\end{equation*}
differentiating (2.8-1) with respect to $t$ again yields
\begin{equation*}
\left.
\begin{array}{ll}
 G'' (t) &\displaystyle =-4Im\int_{\mathbb{R}^N}\left[\dfrac{a_2}{2m_1}x\left(\phi\nabla\overline{\phi}\right)_t
+\dfrac{a_1}{2m_2}x\left(\psi\nabla\overline{\psi}\right)_t\right]dx
\\[0.5cm]
&\displaystyle =-2Im\int_{\mathbb{R}^N}\left[\dfrac{a_2}{m_1}x\left(\phi_t\nabla\overline{\phi}
+\phi\nabla\overline{\phi}_t\right)+\dfrac{a_1}{m_2}x\left(\psi_t\nabla\overline{\psi}
+\psi\nabla\overline{\psi}_t\right)\right]dx.
\end{array}
\right.\eqno(2.8-3)
\end{equation*}
Simple calculations give
\begin{equation*}
\left.
\begin{array}{ll}
&\displaystyle Im\int_{\mathbb{R}^N}\left(x\phi_t\nabla\overline{\phi}+x\phi\nabla
\overline{\phi}_t\right)dx
\\[0.5cm]
 &\displaystyle\quad=Im\int_{\mathbb{R}^N}x\left(\dfrac{i}{2m_1}\Delta\phi
 -ia_1\psi\overline{\phi}\right)\nabla\overline{\phi}dx
+Im\int_{\mathbb{R}^N}x\phi\left(-\dfrac{i}{2m_1}\nabla\Delta\overline{\phi}
+ia_1\nabla(\overline{\psi}\phi)\right)dx
\\[0.5cm]
&\displaystyle\quad=Re\int_{\mathbb{R}^N}\left(\dfrac{x}{2m_1}\Delta\phi\nabla
\overline{\phi}-a_1x\psi\overline{\phi}\nabla\overline{\phi}\right)dx
+Re\int_{\mathbb{R}^N}\left(-\dfrac{1}{2m_1}x\phi\nabla\Delta\overline{\phi}
+a_1x\phi\nabla\left(\overline{\psi}\phi\right)\right)dx,
\end{array}
\right.\eqno(2.8-4)
\end{equation*}
\\
\begin{equation*}
\left.
\begin{array}{ll}
&\displaystyle Im\int_{\mathbb{R}^N}\left(x\psi_t\nabla\overline{\psi}+x\psi(\nabla\overline{\psi})
_t\right)dx
\\[0.5cm]
 &\quad=\displaystyle Im\int_{\mathbb{R}^N}x\left(\dfrac{i}{2m_2}\Delta\psi-ia_2\phi^2\right)\nabla\overline{\psi}dx
+Im\int_{\mathbb{R}^N}x\psi\left(-\dfrac{i}{2m_2}\nabla\Delta
\overline{\psi}+ia_2\nabla\overline{\phi}^2\right)dx
\\[0.5cm]
&\displaystyle\quad=\int_{\mathbb{R}^N}\left(\dfrac{x}{2m_2}\Delta
\psi\nabla\overline{\psi}-a_2x\phi^2\nabla\overline{\psi}\right)dx
+Re\int_{\mathbb{R}^N}\left(-\dfrac{x}{2m_2}\psi\nabla\Delta\overline{\psi}
+a_2x\psi\nabla\overline{\phi}^2\right)dx.
\end{array}
\right.\eqno(2.8-5)
\end{equation*}
\\
Combining (2.8-2), (2.8-3), (2.8-4) and (2.8-5)  together yields\\
\begin{equation*}
\left.
\begin{array}{ll}
G''(t)&\displaystyle=\dfrac{2a_2}{m_1^2}\int_{\mathbb{R}^N}|\nabla\phi|^2dx+
\dfrac{2a_1}{m_2^2}\int_{\mathbb{R}^N}|\nabla\psi|^2dx
\\[0.5cm]
&\quad \displaystyle+2a_1a_2N\left(\dfrac{1}{m_1}-\dfrac{1}{m_2}\right)
Re\int_{\mathbb{R}^N}\overline{\psi}\phi^2dx+2a_1a_2\left(\dfrac{1}{m_1}
-\dfrac{2}{m_2}\right)Re\int_{\mathbb{R}^N}x\nabla\phi^2\overline{\psi}dx.
\end{array}
\right.
\end{equation*}
This completes the proof of  Lemma \ref{Lemma 2.2}.
\end{proof}
We now first list a useful lemma concerning the uniform decay at infinity of certain radial functions in Strass \cite{Strauss1}.
\begin{lem}\label{Lemma 2.3}\cite{Strauss1}( Radial lemma)\quad  Let $N\geq2$. If $u\in H^1(\mathbb{R}^N)$ is a radially symmetric function, then
$$\sup\limits_{x\in \mathbb{R}^N}|x|^{\frac{N-1}{2}}|u(x)|\leq c||u||^{\frac{1}{2}}_{L^2(\mathbb{R}^N)}||\nabla u||^{\frac{1}{2}}_{L^2(\mathbb{R}^N)}\leq c||u||_{H^1(\mathbb{R}^N)}.$$
If, in addition, $u(x)$ is a non-increasing function of $|x|$, then
$$\sup\limits_{x\in \mathbb{R}^N}|x|^{\frac{N}{2}}|u(x)|\leq c||u||_{L^2(\mathbb{R}^N)}.$$
\end{lem}
\indent The following inequality is frequently used throughout this paper.\\
\begin{lem}\label{Lemma 2.4}(Gagliardo-Nirenberg inequality)\quad Let $1\leq p,q,r\leq \infty$, and let $j,m$ be two integers with $0\leq j<m$. If
$$j-\frac{N}{p}=\left(1-\theta\right)\left(-\frac{N}{q}\right)
+\theta\left(m-\frac{N}{r}\right),$$
for some $\theta\in[0,1]$, then there exists $c=c\left(N,m,j,\theta,q,r\right)$ such that
$$\left\|D^ju\right\|_{L^p(\mathbb{R}^N)}\leq c\left\|u\right\|^{1-\theta}_{L^q(\mathbb{R}^N)}\left\|D^mu\right\|^{\theta}_{L^r(\mathbb{R}^N)}$$
for every $u\in \mathcal{D}(\mathbb{R}^N)$.
\end{lem}
We recall below some well-known inequalities and Soboleve embedding results \cite{Caze1}.
\begin{lem}\label{Lemma 2.5}(Poincar\'{e}'s inequality)\quad Assume $|\Omega|<\infty$ (or $\Omega$ is bounded in one direction) and $1\leq p\leq\infty$. Then there exists a constant $c$ such that
$$\left\|u\right\|_{L^p(\Omega)}\leq c\left\|\nabla u\right\|_{L^p(\Omega)}$$
for every $u\in W_0^{1,p}(\Omega)$. In
particular, $||\nabla u||_{L^p(\Omega)}$ is an equivalent norm to $||u|| W^{1,p}(\Omega)$ on $W_0^{1,p}(\Omega)$.
\end{lem}
\begin{lem}\label{Lemma 2.6}(Sobolev's embedding theorem)\quad
If $\Omega$ has a Lipschitz continuous boundary, then the following properties hold:
\\[0.3cm]
(i) If $p>N$, then $W^{1,p}(\Omega)\hookrightarrow L^{\infty}(\Omega)$ .
\\[0.3cm]
If $\Omega$ has a uniformly Lipschitz continuous boundary, then
\\[0.3cm]
(ii) If $p>N$, then $W^{1,p}(\Omega)\hookrightarrow C^{0,\alpha}(\overline{\Omega})$, where $\alpha=\dfrac{p-N}{p}$.
\end{lem}

\begin{lem}\label{Lemma 2.7}(Rellich's compactness theorem)\quad
If $\Omega$ is bounded and has a Lipschitz continuous boundary, then the following properties hold:
\\[0.3cm]
(i) If $p>N$, then the embedding $W^{1,p}(\Omega)\hookrightarrow L^{\infty}(\Omega)$ is compact.
\\[0.3cm]
Let in addition   $\Omega$ have uniformly Lipschitz continuous boundary.
 \\[0.3cm]
(ii) If $p>N$, then the embedding $W^{1,p}(\Omega)\hookrightarrow C^{0,\lambda}(\overline{\Omega})$ is compact for every  $\lambda\in\left(0,\dfrac{p-N}{p}\right)$.
\\[0.3cm]
Furthermore, suppose thay $\Omega$ satisfies the strong local Lipschitz condition.
If $mp>N>(m-1)p$, then\\
$$W^{j+m,p}(\Omega)\hookrightarrow C^{j,\lambda}(\overline{\Omega})$$
for $0<\lambda\leq m-\dfrac{N}{p}$.
\end{lem}

 \begin{lem}\label{Lemma 2.8}(Radial compact lemma)\quad The injection $H_r^1(\mathbb{R}^N)\hookrightarrow L_r^p(\mathbb{R}^N)$ is compact for $2<p<\dfrac{2N}{N-2}$, where $H_r^1(\mathbb{R}^N)$ is the space of all radial functions on $H^1(\mathbb{R}^N)$ and $L_r^p(\mathbb{R}^N)$ the space of all radial functions on $L^p(\mathbb{R}^N)$.
\end{lem}
\section{Finite Time Blowup}
\indent~~ In this section, we show that, under suitable assumptions on initial data, solution of the nonlinear Schr\"{o}dinger equations \eqref{1.1} with this initial data blowups in finite time. We adopt here essentially a convexity analysis method (see Glassey \cite{Glassey1}), which is based on the calculation of the variance
$$\int_{\mathbb{R}^N}|x|^2\left(a_{2}|\phi(t,x)|^2+a_{1}|\psi(t,x)|^2\right)dx.$$
As mentioned at the beginning of section 2, the initial-value problem \eqref{1.1}
-\eqref{2.1} is locally well-posed in $H^1(\mathbb{R}^N)\times H^1(\mathbb{R}^N)$ , and keeps conservation of mass, energy and momentum. Furthermore, the blowup result is based on attaining a suitable virial-type identity. \\
\indent We then clam:
\begin{thm}\label{Theorem 3.1}
~~ Let $4\leq N< 6$; $m_2=2m_1>0, a_1>0, a_2>0$; and let
 $$(\phi,\psi)\in C\left([0,T);H^1(\mathbb{R}^N)\times H^1(\mathbb{R}^N)\right)$$ be the solution to the Cauchy problem \eqref{1.1}-\eqref{2.1}. If in addition,  $\left(|x|\phi_0,|x|\psi_0\right)\in L^2(\mathbb{R}^N)\times L^2(\mathbb{R}^N)$ and either
\\[0.3cm]
(c-1) $E(\phi_0,\psi_0)<0$;
\\[0.3cm]
or
\\[0.3cm]
(c-2) $E(\phi_0,\psi_0)=0$ and $\displaystyle Im\int_{\mathbb{R}^N}\left(\dfrac{a_2}{m_1}x\phi_0\nabla\overline{\phi}_0
+\dfrac{a_1}{m_2}x\psi_0\nabla\overline{\psi}_0\right)dx>0;$
\\[0.3cm]
or
\\[0.3cm]
(c-3) $E(\phi_0,\psi_0)>0$ and
$$Im\int_{\mathbb{R}^N}\left(\dfrac{a_2}{m_1}x\phi_0\nabla\overline{\phi}_0
+\dfrac{a_1}{m_2}x\psi_0\nabla\overline{\psi}_0\right)dx\geq
\left[\dfrac{N}{2m_1}E(\phi_0,\psi_0)\int_{\mathbb{R}^N}|x|^2
\left(a_2|\phi_0|^2+a_1|\psi_0|^2\right)dx\right]^{\frac{1}{2}};$$
\\
then there holds $T<+\infty$ and
 $$\lim\limits_{t\rightarrow T}\left(a_2^{\frac{1}{2}}\|\phi\|_{H^1(\mathbb{R}^N)}
 +a_1^{\frac{1}{2}}\|\psi\|_{H^1(\mathbb{R}^N)}\right)=+\infty.$$
\end{thm}
\begin{proof}
We prove it by contradiction. Suppose that the maximal existence time $T$ of the solution $(\phi,\psi)$ to the Cauchy problem \eqref{1.1}-\eqref{2.1} is infinity. Let
\begin{equation}\label{3.1}
\left.
\begin{array}{ll}
\displaystyle G(t)=\int_{\mathbb{R}^N}|x|^2\left(a_2|\phi|^2+a_1|\psi|^2\right)dx.
\end{array}
\right.
\end{equation}
It follows from \eqref{2.2}, \eqref{2.3},\eqref{2.7} and \eqref{2.8} that
\begin{equation}\label{3.2}
\left.
\begin{array}{ll}
\displaystyle G'(t) =-4Im\int_{\mathbb{R}^N}\left(\dfrac{a_2}{2m_1}x\phi\nabla\overline{\phi}
+\dfrac{a_1}{2m_2}x\psi\nabla\overline{\psi}\right)dx,
\end{array}
\right.
\end{equation}
and
\begin{equation}\label{3.3}
\left.
\begin{array}{ll}
G''(t)&\displaystyle=\dfrac{2a_2}{m_1^2}\int_{\mathbb{R}^N}|\nabla\phi|^2dx
+\dfrac{2a_1}{m_2^2}\int_{\mathbb{R}^N}|\nabla\psi|^2dx
+2a_1a_2N\left(\dfrac{1}{m_1}-\dfrac{1}{m_2}\right)Re\int_{\mathbb{R}^N}
\overline{\psi}\phi^2dx
\\[0.5cm]
&\displaystyle=2N\left(\dfrac{1}{m_1}-\dfrac{1}{m_2}\right)E(\phi_0,\psi_0)
\\[0.5cm]
&\displaystyle\quad-2N\left(\dfrac{1}{m_1}-\dfrac{1}{m_2}\right)
\left(\dfrac{a_2}{2m_1}\int_{\mathbb{R}^N}
|\nabla\phi|^2dx+\dfrac{a_1}{4m_2}\int_{\mathbb{R}^N}
|\nabla\psi|^2dx\right)
\\[0.5cm]
&\displaystyle\quad+\dfrac{2a_2}{m_1^2}\int_{\mathbb{R}^N}|\nabla\phi|^2dx
+\dfrac{2a_1}{m_2^2}\int_{\mathbb{R}^N}|\nabla\psi|^2dx
\\[0.5cm]
&\displaystyle=\dfrac{N}{m_1}E(\phi_0,\psi_0)+\dfrac{a_2}{m_1^2}
\left(2-\dfrac{N}{2}\right)\int_{\mathbb{R}^N}|\nabla\phi|^2dx
+\dfrac{a_1}{4m_1^2}\left(2-\frac{N}{2}\right)\int_{\mathbb{R}^N}|\nabla\psi|^2dx
\\[0.5cm]
&\displaystyle\leq \dfrac{N}{m_1}E(\phi_0,\psi_0),
\end{array}
\right.
\end{equation}
where the condition $m_2=2m_1>0$ ensures the second to last equality valid, while the last one holds true due to the condition $ N\geq 4$. On the other hand, classical discussion then yields
\begin{equation}\label{3.4}
G(t)=G(0)+G'(0)t+\int_0^t(t-s)G''(s)dx, \quad 0\leq t<+\infty.
\end{equation}
From \eqref{3.3} it follows that
\begin{equation}\label{3.5}
G(t)\leq G(0)+G'(0)t+\dfrac{N}{2m_1}E(\phi_0,\psi_0)t^2, \quad 0\leq t<\infty.
\end{equation}
Noting that $G(t)$ is a nonnegative function and
\begin{equation}\label{3.6}
G(0)=\int_{\mathbb{R}^N}|x|^2\left(a_2|\phi_0|^2+a_1|\psi_0|^2\right)dx\geq0,
\end{equation}
\eqref{3.2} thus yields
\begin{equation}\label{3.7}
G'(0)=-2Im\int_{\mathbb{R}^N}\left(\dfrac{a_2}{m_1}x\phi_0\nabla\overline{\phi_0}
+\dfrac{a_1}{m_2}x\psi_0\nabla\overline{\psi_0}\right)dx.
\end{equation}
\\
Hence, under the assumptions (c-1) or (c-2) or (c-3), \eqref{3.5} implies that there exists $T^*<+\infty$ such that
\begin{equation}\label{3.8}
\lim\limits_{t\rightarrow T^*}G(t)=\lim\limits_{t\rightarrow T^*}\int_{\mathbb{R}^N}|x|^2\left(a_2|\phi|^2+a_1|\psi|^2\right)dx=0.
\end{equation}
By H\"{o}lder's inequality and Schwarz inequality, we have
\begin{equation}\label{3.9}
\left.
\begin{array}{ll}
&\displaystyle~~a_2\int_{\mathbb{R}^N}|\phi|^2dx+a_1\int_{\mathbb{R}^N}|\psi|^2dx
\qquad\qquad\qquad
\\[0.5cm]
&\quad\displaystyle=-\dfrac{2}{N}Re\int_{\mathbb{R}^N}a_2x\phi\nabla\overline{\phi}dx
- \dfrac{2}{N}Re\int_{\mathbb{R}^N}a_1x\psi\nabla\overline{\psi}dx
 \\[0.5cm]
&\quad\displaystyle\leq \dfrac{2}{N}\left\|a_2^{\frac{1}{2}}x\phi\right\|_{L^2(\mathbb{R}^N)}
\left\|a_2^{\frac{1}{2}}\nabla\phi\right\|_{L^2(\mathbb{R}^N)}
+\dfrac{2}{N}\left\|a_1^{\frac{1}{2}}x\psi\right\|_{L^2(\mathbb{R}^N)}
\left\|a_1^{\frac{1}{2}}
\nabla\psi\right\|_{L^2(\mathbb{R}^N)}
\\[0.5cm]
&\quad\displaystyle\leq\dfrac{4}{N}\left(\left\|a_2^{\frac{1}{2}}x\phi\right\|
^2_{L^2(\mathbb{R}^N)}+\left\|a_1^{\frac{1}{2}}x\psi\right\|^2_{L^2(\mathbb{R}^N)}\right)^{\frac{1}{2}}
\left(\left\|a_2^{\frac{1}{2}}\nabla\phi\right\|^2_{L^2(\mathbb{R}^N)}+\left\|a_1^{\frac{1}{2}}
\nabla\psi\right\|^2_{L^2(\mathbb{R}^N)}\right)^{\frac{1}{2}}
\\[0.5cm]
&\quad\displaystyle\leq\dfrac{4}{N}\left(\int_{\mathbb{R}^N}a_2|x|^2|\phi|^2dx
+\int_{\mathbb{R}^N}a_1|x|^2|\psi|^2dx\right)^{\frac{1}{2}}
\left(\int_{\mathbb{R}^N}a_2|\nabla\phi|^2dx+\int_{\mathbb{R}^N}
a_1|\nabla\psi|^2dx\right)^{\frac{1}{2}}.
\end{array}
\right.
\end{equation}
Therefore, as $t\rightarrow T^*$, \eqref{3.9} together with \eqref{3.8} implies that
$$a_2\int_{\mathbb{R}^N}|\phi|^2dx+a_1\int_{\mathbb{R}^N}|\psi|^2dx\leq 0,$$
which is a contradiction from the conservation of mass \eqref{2.2}.\\
\indent This finishes the proof of  Theorem \ref{Theorem 3.1}.
\end{proof}
\section{The Existence of Standing Waves associated with ground state}
\indent \quad Under the assumption $\omega_2=2\omega_1>0$, $m_2=2m_1$ and $ 4<N< 6$, if $\left(u(x),v(x)\right)\in H^1(\mathbb{R}^N)\times H^1(\mathbb{R}^N)\backslash \{(0,0)\}$ is the ground state solution to \eqref{1.5}, then $\left(\phi(t,x),\psi(t,x)\right)\equiv\left(e^{i\omega_1t}u(x),
e^{i\omega_2t}v(x)\right)$ is a standing wave solution of \eqref{1.1}. Hence, it is sufficient to explore the existence of the ground state solution of \eqref{1.5}.\\
\indent Rewriting equations \eqref{1.5} as the following:
\begin{equation}\label{4.1}
\left\{
\begin{array}{ll}
-\dfrac{1}{2m_1}\Delta u(x)=-\omega_1 u(x)-a_1u(x)v(x)\triangleq g_1(u(x),v(x)),\quad(4.1-a)
\\[0.5cm]
-\dfrac{1}{2m_2}\Delta v(x)=-\omega_2 v(x)-a_2u^2(x)\triangleq g_2(u(x),v(x)),\quad(4.1-b)
\end{array}
\right.
\end{equation}
concerning the existence and exponential decay at infinity of the ground state solution for equations \eqref{4.1}, we then claim
\begin{thm}\label{Theorem 4.1}
Let $\omega_2=2\omega_1>0$, $m_2=2m_1$, $4<N<6$ and \eqref{1.10} hold true. There exists $\left(\xi(x),-\eta(x)\right)\in\mathcal{M}$ with $\xi(x)>0, ~~\eta(x)>0$ such that
\\[0.3cm]
\indent (1) \quad $\mathcal{S}\left(\xi(x),-\eta(x)\right)=K=\inf\limits_{u,v\in\mathcal{M}}\mathcal{S}(u,v)$;
\\[0.3cm]
\indent (2)\quad $\left(\xi(x),-\eta(x)\right)$ is a ground state solution of \eqref{4.1};
\\[0.3cm]
\indent (3)\quad $\xi(x)$ and $\eta(x)$ are functions of $|x|$ alone and have exponential decays at infinity.
\\[0.3cm]
Here, $\mathcal{S}$ and $\mathcal{M}$ are defined by \eqref{1.6} and \eqref{1.8}, respectively.
\end{thm}
Before proving Theorem \ref{Theorem 4.1}, we first give a key conclusion.
\begin{lem}\label{Lemma 4.2}
Assuming that \eqref{1.10} holds true, let $\omega_2=2\omega_1>0$, $m_2=2m_1>0$, $4<N<6$, and let
$(u,v)\in H^1(\mathbb{R}^N)\times H^1(\mathbb{R}^N)$ be the solution of \eqref{4.1}. Then the functions
$$g_1(u,v)=-\omega_1 u(x)-a_1u(x)v(x) ~~\mbox{and}~~ g_2(u,v)=-\omega_2 v(x)-a_2u^2(x)$$
satisfy the following conditions.\\
\indent For $L=\dfrac{N+2}{N-2}$, there hold
 $$1-L=\dfrac{-4}{N-2}, ~~2< L<3,~~ 0< L-2<1,\eqno (S-1)$$
$$-\infty<\liminf\limits_{(u,v)\rightarrow(0^+,0^-)}\dfrac{g_1(u,v)}{u}\leq \limsup\limits_{(u,v)\rightarrow(0^+,0^-)}\dfrac{g_1(u,v)}{u}=-\omega_1<0,\eqno (S-2)$$
$$-\infty<\liminf\limits_{(u,v)\rightarrow(0^+,0^-)}\dfrac{g_2(u,v)}{v}\leq \limsup\limits_{(u,v)\rightarrow(0^+,0^-)}\dfrac{g_2(u,v)}{v}=-\omega_2<0,\eqno  (S-3)$$
$$-\infty\leq\limsup\limits_{{(u,v)\rightarrow(+\infty,-\infty)}}
\dfrac{g_1(u,v)}{u^L}\leq0,\eqno (S-4)$$
$$-\infty\leq\limsup\limits_{{(u,v)\rightarrow(+\infty,-\infty)}}
\dfrac{g_2(u,v)}{v^L}\leq0.\eqno (S-5)$$
In particular, $g_1(u,v)$ and  $g_2(u,v)$ satisfy two stronger conditions for $L=\dfrac{N+2}{N-2}$:
$$\lim\limits_{(u,v)\rightarrow(\pm\infty,\pm\infty)}\dfrac{|g_1(u,v)|}{|u|^L}=0,\eqno (S-4)^*$$
$$\lim\limits_{(u,v)\rightarrow(\pm\infty,\pm\infty)}\dfrac{|g_2(u,v)|}{|v|^L}=0.\eqno (S-5)^*$$
Furthermore, there exist $u^*>0$ and $v^*<0$ such that
$$G(u^*,v^*)=-a_1a_2{u^*}^2v^*-a_2\omega_1{u^*}^2-\dfrac{a_1\omega_2{v^*}^2}{2}>0,\eqno  (S-6)$$
where $$\dfrac{\partial G(u^*,v^*)}{\partial u^*}=-2a_2\omega_1u^*-2a_1a_2u^*v^*,\eqno (S-7)$$
and $$\dfrac{\partial G(u^*,v^*)}{\partial v^*}=-a_1\omega_2v^*-a_1a_2{u^*}^2.\eqno  (S-8) $$
\end{lem}
\begin{proof}
Since $4< N<6$, for $L=\dfrac{N+2}{N-2}$, direct calculation leads to (S-1). Noticing that\\
$$g_1(u,v)=-\omega_1u-a_1uv,\quad g_2(u,v)=-\omega_2v-a_2u^2,$$
\\ by \eqref{4.1}, and using L'Hospital's rule, we have
\begin{equation*}\label{S-9}
\left.
\begin{array}{ll}
-\infty&\displaystyle<\liminf\limits_{(u,v)\rightarrow(0^+,0^-)}\dfrac{g_1(u,v)}{u}  \leq \limsup\limits_{(u,v)\rightarrow(0^+,0^-)}\dfrac{g_1(u,v)}{u}
\\[0.5cm]
&\displaystyle=\limsup\limits_{(u,v)\rightarrow(0^+,0^-)}\dfrac{-\omega_1 u-a_1uv}{u}=\limsup\limits_{(u,v)\rightarrow(0^+,0^-)}(-\omega_1-a_1v)=-\omega_1<0,
\end{array}
\right.\eqno(S-9)
\end{equation*}\\
\begin{equation*}\label{S-10}
\left.
\begin{array}{ll}
-\infty&<\liminf\limits_{(u,v)\rightarrow(0^+,0^-)}\dfrac{g_2(u,v)}{v}
 \leq \limsup\limits_{(u,v)\rightarrow(0^+,0^-)}\dfrac{g_2(u,v)}{v}
\\[0.5cm]
&\displaystyle
 =\limsup\limits_{(u,v)\rightarrow(0^+,0^-)}\dfrac{-\omega_2 v-a_2u^2}{v}=\limsup\limits_{(u,v)\rightarrow(0^+,0^-)}\left(-\omega_2
 -a_2\dfrac{u^2}{v}\right)
 \\[0.5cm]
&\displaystyle
 =-\omega_2-a_2\limsup\limits_{(u,v)\rightarrow(0^+,0^-)}\dfrac{u^2}{v}=-\omega_2<0,
\end{array}
\right.\eqno(S-10)
\end{equation*}\\
\begin{equation*}\label{S-11}
\left.
\begin{array}{ll}
-\infty&\displaystyle\leq\liminf\limits_{(u,v)\rightarrow(+\infty,-\infty)}
\dfrac{g_1(u,v)}{u^L} \leq \limsup\limits_{(u,v)\rightarrow(+\infty,-\infty)}\dfrac{g_1(u,v)}{u^L}
  \\[0.5cm]
&\displaystyle=\limsup\limits_{(u,v)\rightarrow(+\infty,-\infty)}\dfrac{-\omega_1 u-a_1uv}{u^L}=\limsup\limits_{(u,v)\rightarrow(+\infty,-\infty)}
\left(\dfrac{-\omega_1}{u^{L-1}}-a_1\dfrac{v}{u^{L-1}}\right)
 \\[0.5cm]
&\displaystyle=\limsup\limits_{(u,v)\rightarrow(+\infty,-\infty)}\left
(-a_1\dfrac{v}{u^{L-1}}\right)
=\limsup\limits_{(u,v)\rightarrow(+\infty,-\infty)}\left(\dfrac{-a_1}
{(L-1)u^{L-2}}\right)=0,
\end{array}
\right.\eqno(S-11)
\end{equation*}\\
\begin{equation*}\label{S-12}
\left.
\begin{array}{ll}
-\infty&\displaystyle\leq\liminf\limits_{(u,v)\rightarrow(+\infty,-\infty)}
\dfrac{g_2(u,v)}{v^L}
 \leq \limsup\limits_{(u,v)\rightarrow(+\infty,-\infty)}\dfrac{g_2(u,v)}{v^L}
\\[0.5cm]
&\displaystyle=\limsup\limits_{(u,v)\rightarrow(+\infty,-\infty)}\left
(\dfrac{-\omega_2 }{v^{L-1}}-a_2\dfrac{u^2 }{v^L}\right)
=\limsup\limits_{(u,v)\rightarrow(+\infty,-\infty)}
\left(\dfrac{-a_2}{L(L-1)v^{L-2}}\right)=0,
\end{array}
\right.\eqno(S-12)
\end{equation*}\\
\begin{equation*}\label{S-13}
\left.
\begin{array}{ll}
\displaystyle\lim\limits_{(u,v)\rightarrow(\pm\infty,\pm\infty)}\dfrac{|g_1(u,v)|}{|u|^L}
&\displaystyle=\lim\limits_{(u,v)\rightarrow(\pm\infty,\pm\infty)}\dfrac{|-\omega_1u-a_1uv|}{|u|^L}
\\[0.5cm]
&\displaystyle\leq\lim\limits_{(u,v)\rightarrow(\pm\infty,\pm\infty)}\left
(\frac{\omega_1|u|}{|u|^L}+\frac{a_1|u||v|}{|u|^L}\right)
 =0,
\end{array}
\right.\eqno(S-13)
\end{equation*}\\
\begin{equation*}\label{S-14}
\left.
\begin{array}{ll}
\displaystyle\lim\limits_{(u,v)\rightarrow(\pm\infty,\pm\infty)}\dfrac{|g_2(u,v)|}{|v|^L}
&\displaystyle=\lim\limits_{(u,v)\rightarrow(\pm\infty,\pm\infty)}\dfrac{|-\omega_2v-a_2u^2|}{|v|^L}
\\[0.5cm]
&\displaystyle\leq\lim\limits_{(u,v)\rightarrow(\pm\infty,\pm\infty)}
\left(\dfrac{\omega_2|v|}{|v|^L}+\dfrac{a_2|u|^2}{|v|^L}\right)
 =0.
\end{array}
\right.\eqno(S-14)
\end{equation*}\\
Combining $(S-9)$ with $(S-10),~ (S-11), ~(S-12), ~( S-13)$ and $(S-14)$ yields $(S-1),~ (S-2), ~(S-3), ~(S-4),~(S-5),~ (S-4)^*$ and $(S-5)^*$.\\
\indent We then verify estimate $(S-6)$. Using $(S-7)$, one obtains\\
\begin{equation*}\label{S-15}
\left.
\begin{array}{ll}
\displaystyle G(u^*,v^*)&\displaystyle=\int_0^{u^*}\left(-2a_2\omega_1s-2a_1a_2sv^*\right)ds+\Phi(v^*)
\\[0.5cm]
&\displaystyle=-a_2\omega_1{u^*}^2-a_1a_2{u^*}^2v^*+\Phi(v^*).
\end{array}
\right.\eqno(S-15)
\end{equation*}\\
Since $(u^*,v^*)$ satisfies $(S-8)$, differentiating $(S-15)$ with respect to $v^*$ yields
$$-a_1a_2{u^*}^2+\dfrac{d\Phi(v^*)}{dv^*}=-a_1\omega_2{v^*}-a_1a_2{u^*}^2,$$
that is,
\begin{equation*}
\left.
\begin{array}{ll}
\displaystyle\dfrac{d\Phi(v^*)}{dv^*}=-a_1\omega_2{v^*},
\end{array}
\right.
\end{equation*}
then
\begin{equation*}\label{S-16}
\left.
\begin{array}{ll}
\displaystyle\Phi(v^*)=-\dfrac{a_1}{2}\omega_2{v^*}^2.
\end{array}
\right.\eqno(S-16)
\end{equation*}
This together with $(S-16)$ leads to
\begin{equation*}\label{S-17}
\left.
\begin{array}{ll}
\displaystyle G(u^*,v^*)=-a_2\omega_1{u^*}^2-\dfrac{a_1}{2}\omega_2{v^*}^2-a_1a_2{u^*}^2v^*.
\end{array}
\right.\eqno(S-17)
\end{equation*}
Since $\left\{(u,v):G(u,v)=0\right\}$ is an one-dimensional set and $$G(u,-u)=-a_2\omega_1u^2-\dfrac{a_1}{2}\omega_2u^2+a_1a_2u^3,$$
we have $G(u,-u)>0$ if $u\geq u_0$ for $u_0$ large enough. Thus we can find $u^*>0$ and $v^*<0$ such that
$$G(u^*,v^*)=-a_2\omega_1{u^*}^2-\frac{a_1}{2}\omega_2{v^*}^2-a_1a_2{u^*}^2v^*>0.$$
This completes the proof of Lemma \ref{Lemma 4.2}.
\end{proof}
\indent We now begin to show Theorem \ref{Theorem 4.1}, which will be divided into three subsections.\\
\begin{proof}
\indent 4.1. Existence of solutions to the minimizing problem \eqref{1.9} \Big(~proof of (1) in Theorem \ref{Theorem 4.1}~\Big).\\
\indent 4.2. Existence of ground state solution of  \eqref{4.1} \Big(~proof of (2) in Theorem \ref{Theorem 4.1}~\Big).\\
\indent 4.3. Exponential decay of the ground state solution of \eqref{4.1} \Big(~proof of (3) in Theorem \ref{Theorem 4.1}~\Big).
\subsection{Existence of solutions to the minimizing problem \eqref{1.9}
\\[0.3cm]  \indent ----proof of (1) in Theorem \ref{Theorem 4.1} }
\indent\quad In order to prove (1) in Theorem \ref{Theorem 4.1}, it is sufficient to show the following result which provides a variational characterization of the ground states to \eqref{4.1}.
\begin{prop}\label{Proposition 4.3}
Let $\omega_2=2\omega_1$ and  $4< N<6$, assuming that \eqref{1.10} holds true. Then any solution $(u,v)$ of \eqref{4.1} belongs to $\mathcal{M}$.  In addition, there exists $(\xi, -\eta)\in\mathcal{M}$ with $\xi>0, ~~\eta>0$ such that

\begin{equation}\label{4.2}
\left.
\begin{array}{ll}
\mathcal{S}(\xi,-\eta)=\min\limits_{(u,v)\in\mathcal{M}}\mathcal{S}(u,v).
\end{array}
\right.
\end{equation}
Here,  $\mathcal{S}$ and $\mathcal{M}$ are defined by \eqref{1.6} and \eqref{1.8},
respectively.
\\[0.3cm]
{\bf Remark 4.3.1.} Theorem \ref{Theorem 4.1} indicates that any ground state solution of \eqref{4.1} is a solution of the minimization problem \eqref{4.2}.
\end{prop}
\indent Before proving Proposition {\ref{Proposition 4.3}} , we recall here the basic properties of Schwarz symmetrization. We first mention the definition of the Schwarz spherical rearrangement (or symmetrization) of a function (see Berestycki-Lions\cite{Bere-Lions1}).
\begin{defn}\label{Definition 4.4}\cite{Bere-Lions1}(Schwarz symmetrization )~~ Let $f\in L^1(\mathbb{R}^N)$ be a nonnegative function, then $f^*$, the Schwarz symmetrized function of $f$, is the unique spherically symmetric, non-increasing (in $r=|x|$), measurable function such that for all $\alpha>0$,
$m\left\{x\in\mathbb{R}^N: f^*\geq\alpha\right\}=m\left\{x\in\mathbb{R}^N: |f|\geq\alpha\right\}$, where $m$ is the Lebesgue measure.
\end{defn}
We next refer to Berstycki and Lions \cite{Bere-Lions1}, Appendix AIII for the main properties of the Schwarz symmetrization.
\begin{lem}\label{Lemma 4.5}\cite{Bere-Lions1}(Basic properties of Schwarz symmetrization)
~Let $f^*$ and $g^*$ be the Schwarz symmetrization of functions $ |f| $ and $ |g| $, respectively. Then there hold:
\\[0.3cm]
 (1) For every continuous function $F$ such that $F(f)$ is integrable, then $$\int_{\mathbb{R}^N}F(f)dx=\int_{\mathbb{R}^N}F(f^*)dx.$$
\\
 (2)~ Riesz inequality:   Let $f^*$ and $g^*$ be in $L^2(\mathbb{R}^N)$, then $$\int_{\mathbb{R}^N}f(x)g(x)dx\leq\int_{\mathbb{R}^N}f^*(x)g^*(x)dx.$$
\\
 (3)~ Let $f$ , $g$ be in $L^2(\mathbb{R}^N)$, then $||f^*-g^*||_{L^2(\mathbb{R}^N)}\leq||f-g||_{L^2(\mathbb{R}^N)}$.
\\[0.3cm]
 (4)~ $\displaystyle\int_{\mathbb{R}^N}|f^*|^pdx=\int_{\mathbb{R}^N}|f|^pdx$ for all~$1\leq p<\infty$ such that $f\in L^p(\mathbb{R}^N)$.
\\[0.3cm]
 (5)~ Let $f$ be in $\mathcal{D}^{1,2}(\mathbb{R}^N)$ if $N\geq 3$~ \Big(respectively, in $H^1(\mathbb{R}^N)$ for any $N$\Big), then $f^*$  belongs to $\mathcal{D}^{1,2}(\mathbb{R}^N)$~ \Big(respectively, to $H^1(\mathbb{R}^N)$\Big), and there holds ~$\displaystyle\int_{\mathbb{R}^N}|\nabla f^*|^2dx\leq\int_{\mathbb{R}^N}|\nabla f|^2dx$.
\\[0.3cm]
 (6)~ Let $f_{\lambda}(x)=\lambda^{\frac{N}{2}}f(\lambda x)$,~ then $(f_{\lambda})^*=(f^*)_{\lambda}$.
\end{lem}
From the definitions of $\mathcal{S}(u,v)$, $\mathcal{Q}(u,v)$ and $\mathcal{M}$ formulated by \eqref{1.6},
\eqref{1.7} and \eqref{1.8}, respectively, we claim
\begin{lem}\label{Lemma 4.6}
Let $\omega_2=2\omega_1>0$, $4< N<6$ and \eqref{1.10} hold true. Then $\mathcal{S}(u,v)$ is bounded below on $\mathcal{M}$.
\end{lem}
\begin{proof}
If $(u,v)\in\mathcal{M}$, then from \eqref{1.6}, \eqref{1.7} and \eqref{1.8} it follows that
\begin{equation}\label{4.3}
\left.
\begin{array}{ll}
\mathcal{S}(u,v)&\displaystyle=\dfrac{N-4}{2N}\left(\dfrac{a_2}{m_1}
\int_{\mathbb{R}^N}|\nabla u|^2dx+\dfrac{a_1}{2m_2}\int_{\mathbb{R}^N}|\nabla v|^2dx\right)
\\[0.5cm]
&~~\displaystyle+a_2\omega_1\int_{\mathbb{R}^N}|u|^2dx+\dfrac{a_1}{2}\omega_2\int_{\mathbb{R}^N}|v|^2dx.
\end{array}
\right.
\end{equation}
Note that $\omega_2=2\omega_1>0$ and $N>4$, \eqref{4.3} yields that $\mathcal{S}(u,v)>0$ for any $(u,v)\in\mathcal{M}$.
\end{proof}
We then formulate a technique lemma.
\begin{lem}\label{Lemma 4.7}
Let $\omega_2=2\omega_1>0$, and  $4< N< 6$, assuming that \eqref{1.10}
holds true.  For
$$(u,v)\in H^1(\mathbb{R}^N)\times H^1(\mathbb{R}^N)\setminus\{(0,0)\}\quad and\quad\lambda>0,$$ let
\begin{equation}\label{4.4}
\left.
\begin{array}{ll}
\displaystyle u_{\lambda}(x)=\lambda^{\frac{N}{2}}u(\lambda x),\quad v_{\lambda}(x)=\lambda^{\frac{N}{2}}v(\lambda x).
\end{array}
\right.
\end{equation}
Then there exists a unique $\beta>0$ (depending on $(u,v)$) such that $Q(u_{\beta},v_{\beta})=0$, ~and
\begin{equation}\label{4.5}
\left\{
\begin{array}{ll}
Q(u_{\lambda},v_{\lambda})>0&\quad for \quad \lambda\in(0,\beta),
\\[0.3cm]
Q(u_{\lambda},v_{\lambda})<0&\quad for \quad \lambda\in(\beta,\infty),
\\[0.3cm]
\mathcal{S}(u_{\beta},v_{\beta})\geq\mathcal{S}(u_{\lambda},v_{\lambda})&\quad for \quad any\quad \lambda>0.\\
\end{array}
\right.
\end{equation}
\end{lem}
\begin{proof}
From \eqref{1.6} and \eqref{1.7}, direct computation gives
\begin{equation*}
\left.
\begin{array}{ll}
\mathcal{S}(u_{\lambda},v_{\lambda})&\displaystyle=\dfrac{a_2}{2m_1}\lambda^2
\int_{\mathbb{R}^N}|\nabla u|^2dx+\dfrac{a_1}{4m_2}\lambda^2\int_{\mathbb{R}^N}|\nabla v|^2dx
\\[0.5cm]
&\displaystyle~~+a_2\omega_1\int_{\mathbb{R}^N}|u|^2dx+\dfrac{a_1}{2}\omega_2\int_{\mathbb{R}^N}|v|^2dx
+\lambda^{\frac{N}{2}}a_1a_2\int_{\mathbb{R}^N}vu^2dx,
\end{array}
\right.
\end{equation*}\\
\begin{equation*}
\left.
\begin{array}{ll}
 Q(u_{\lambda},v_{\lambda})
 & \displaystyle=\dfrac{a_2}{m_1}\lambda^2\int_{\mathbb{R}^N}|\nabla u|^2dx+\dfrac{a_1}{2m_2}\lambda^2\int_{\mathbb{R}^N}|\nabla v|^2dx
+\dfrac{N}{2}\lambda^{\frac{N}{2}}a_1a_2\int_{\mathbb{R}^N}vu^2dx
\\[0.5cm]
& \displaystyle=\lambda^{\frac{N}{2}}\left[\lambda^{2-\frac{N}{2}}\left
(\dfrac{a_2}{m_1}\int_{\mathbb{R}^N}|\nabla u|^2dx+\dfrac{a_1}{2m_2}\lambda^2\int_{\mathbb{R}^N}|\nabla v|^2dx\right)+\dfrac{N}{2}a_1a_2\int_{\mathbb{R}^N}vu^2dx\right].
\end{array}
\right.
\end{equation*}
\\
Recalling $(S-6)$ in Lemma \ref{Lemma 4.2}, there exists $\beta>0$ such that $Q(u_{\beta},v_{\beta})=0$, where $\beta$ depends on $(u,v)$ with $v<0$.
In addition, there holds
$$Q(u_{\lambda},v_{\lambda})>0\quad for \quad \lambda\in(0,\beta), \quad Q(u_{\lambda},v_{\lambda})<0\quad for \quad \lambda\in(\beta,+\infty).$$
Note that$$\frac{d}{d\lambda}\mathcal{S}(u_{\lambda},v_{\lambda})=\lambda^{-1}Q(u_{\lambda},v_{\lambda}),$$
and$$Q(u_{\beta},v_{\beta})=0,$$
one knows for any $\lambda>0$, $\mathcal{S}(u_{\beta},v_{\beta})\geq \mathcal{S}(u_{\lambda},v_{\lambda})$.\\
\indent This completes the proof of Lemma \ref{Lemma 4.7}.
\end{proof}
We further give a crucial conclusion.
\begin{lem}\label{Lemma 4.8}
Solutions of \eqref{4.1}
belong to $\mathcal{M}$.
\end{lem}
\begin{proof}
Let $(u,v)\in H^1(\mathbb{R}^N)\times H^1(\mathbb{R}^N)\setminus\{(0,0)\}$ be a solution of \eqref{4.1}. Firstly, multiplying $(4.1-a)$ by $2a_2u$ and $(4.1-b)$ by $a_1v$, then integrating over $\mathbb{R}^N$ with respect to $x$ yields
\begin{equation*}\label{B-1}
\left.
\begin{array}{ll}
\displaystyle\dfrac{a_2}{m_1}\int_{\mathbb{R}^N}|\nabla u|^2dx&\displaystyle+\dfrac{a_1}{2m_2}\int_{\mathbb{R}^N}|\nabla v|^2dx+2a_2\omega_1\int_{\mathbb{R}^N}u^2dx
\\[0.5cm]
 & \displaystyle+a_1\omega_2\int_{\mathbb{R}^N}v^2dx+3a_1a_2\int_{\mathbb{R}^N}vu^2dx=0.
\end{array}
\right.\eqno(B-1)
\end{equation*}
On the other hand, multiplying $(4.1-a)$ by $2a_2x\cdot \nabla u$ and $(4.1-b)$ by $a_1x\cdot \nabla v$, then integrating over $\mathbb{R}^N$, we obtain
\begin{equation*}\label{B-2}
\left.
\begin{array}{ll}
\displaystyle\dfrac{N-2}{Nm_1}a_2\int_{\mathbb{R}^N}|\nabla u|^2dx&\displaystyle+\dfrac{N-2}{2Nm_2} a_1\int_{\mathbb{R}^N}|\nabla v|^2dx + 2a_2\omega_1\int_{\mathbb{R}^N}u^2dx
\\[0.5cm]
&\displaystyle+a_1\omega_2\int_{\mathbb{R}^N}v^2dx+2a_1a_2\int_{\mathbb{R}^N}vu^2dx=0.
\end{array}
\right.\eqno(B-2)
\end{equation*}
Now subtracting (B-2) from (B-1) leads to
\begin{equation*}
\left.
\begin{array}{ll}
\displaystyle\dfrac{2a_2}{Nm_1}\int_{\mathbb{R}^N}|\nabla u|^2dx+\dfrac{a_1}{Nm_2}\int_{\mathbb{R}^N}|\nabla v|^2dx+a_1a_2\int_{\mathbb{R}^N}vu^2dx=0,
\end{array}
\right.
\end{equation*}
which is equivalent to
\begin{equation*}
\left.
\begin{array}{ll}
\displaystyle\dfrac{a_2}{m_1}\int_{\mathbb{R}^N}|\nabla u|^2dx+\dfrac{a_1}{2m_2}\int_{\mathbb{R}^N}|\nabla v|^2dx+\frac{N}{2}a_1a_2\int_{\mathbb{R}^N}vu^2dx=0,
\end{array}
\right.
\end{equation*}
that is, $Q(u,v)=0$.\\
\indent This completes the proof of Lemma \ref{Lemma 4.8}.
\end{proof}

We are now in the position to show Proposition \ref{Proposition 4.3}.\\

{\bf Proof of Proposition \ref{Proposition 4.3} }.
\begin{proof} Let $\left\{(u_n,v_n),n\in\mathbb{N}\right\}\subset\mathcal{M}$ be a minimizing sequence for \eqref{1.9}, that is, $(u_n,v_n)\neq (0,0)$, and as $n\rightarrow+\infty$,
\begin{equation}\label{4.6}
\left.
\begin{array}{ll}
\displaystyle\mathcal{S}(u_n,v_n)\rightarrow
\inf\limits_{(u,v)\in\mathcal{M}}\mathcal{S}(u,v),
\end{array}
\right.
\end{equation}
as well as
\begin{equation*}
\left.
\begin{array}{ll}
Q(u_{n},v_{n})=\displaystyle\dfrac{a_2}{m_1}\int_{\mathbb{R}^N}|\nabla u_{n}|^2dx+\dfrac{a_1}{2m_2}\int_{\mathbb{R}^N}|\nabla v_{n}|^2dx
+\dfrac{N}{2}a_1a_2\int_{\mathbb{R}^N}v_{n} u_{n} ^2dx=0,
\end{array}
\right.\eqno(4.6)^*
\end{equation*}
which implies that $v_{n}$ needs to satisfy $v_{n}<0$. Thus, we can rewritten $(4.6)^*$ as
 \begin{equation*}
\left.
\begin{array}{ll}
Q(u_{n},v_{n})=\displaystyle\dfrac{a_2}{m_1}\int_{\mathbb{R}^N}|\nabla u_{n}|^2dx+\dfrac{a_1}{2m_2}\int_{\mathbb{R}^N}|\nabla v_{n}|^2dx
-\dfrac{N}{2}a_1a_2\int_{\mathbb{R}^N}(-v_{n}) u_{n} ^2dx=0.
\end{array}
\right.\eqno(4.6)^a
\end{equation*}
 According to Lemma 4.4 and Lemma 4.5, let $u_{n}^*$ and $v_{n}^*$ be the Schwarz spherical rearrangement of functions $|u_{n}|$ and $|v_{n}|=-v_{n}$ (with $v_{n}<0$), respectively. Then $Q(u^*_{n},-v^*_{n})$ can be written as
\begin{equation*}
\left.
\begin{array}{ll}
Q(u^{*}_{n},-v^{*}_{n})=\displaystyle\dfrac{a_2}{m_1}\int_{\mathbb{R}^N}|\nabla u^{*}_{n}|^2dx+\dfrac{a_1}{2m_2}\int_{\mathbb{R}^N}|\nabla v^*_{n}|^2dx
-\dfrac{N}{2}a_1a_2\int_{\mathbb{R}^N}v^{*}_{n} {u^{*}_{n}} ^2dx.
\end{array}
\right.\eqno(4.7)^a
\end{equation*}
Referring to Lemma 4.4 and Lemma 4.5 again, we obtain
\begin{equation*}
\left.
\begin{array}{ll}
Q(u^{*}_{n},-v^{*}_{n})&=\displaystyle\dfrac{a_2}{m_1}\int_{\mathbb{R}^N}|\nabla u^{*}_{n}|^2dx+\dfrac{a_1}{2m_2}\int_{\mathbb{R}^N}|\nabla v^{*}_{n}|^2dx
-\dfrac{N}{2}a_1a_2\int_{\mathbb{R}^N} v^{*}_{n} {u^{*}_{n} }^2dx
\\[0.5cm]
&\leq \displaystyle\dfrac{a_2}{m_1}\int_{\mathbb{R}^N}|\nabla u_{n}|^2dx+\dfrac{a_1}{2m_2}\int_{\mathbb{R}^N}|\nabla v_{n}|^2dx
-\dfrac{N}{2}a_1a_2\int_{\mathbb{R}^N}(-v _{n} u _{n} ^2)dx
\\[0.5cm]
&\leq \displaystyle\dfrac{a_2}{m_1}\int_{\mathbb{R}^N}|\nabla u_{n}|^2dx+\dfrac{a_1}{2m_2}\int_{\mathbb{R}^N}|\nabla v_{n}|^2dx
+\dfrac{N}{2}a_1a_2\int_{\mathbb{R}^N}v _{n} u _{n} ^2dx
\\[0.5cm]
&=Q(u_{n},v_{n})=0.
\end{array}
\right.\eqno(4.7)^b
\end{equation*}
As given in Lemma \ref{Lemma 4.7}, put
\begin{equation}\label{4.7}
\left.
\begin{array}{ll}
\displaystyle (u_{n})_{\lambda}=\lambda^{\frac{N}{2}}u_{n}(\lambda x),\quad (v_{n})_{\lambda}=\lambda^{\frac{N}{2}}v_{n}(\lambda x).
\end{array}
\right.
\end{equation}
Then for the minimizing sequence $\left\{(u_n,v_n),n\in\mathbb{N}\right\}\subset\mathcal{M}$, we let
$$\xi_n=(u^*_n)_{\beta_n},\quad\eta_n=(v^*_n)_{\beta_n},$$
where $0<\beta_n\leq 1$ is uniquely determined by
\begin{equation}\label{4.8}
\left.
\begin{array}{ll}
\displaystyle Q(\xi_n,-\eta_n)=Q\left((u^*_n)_{\beta_n},-(v^*_n)_{\beta_n}\right)=0.
\end{array}
\right.
\end{equation}
On the other hand, by $(6)$ of Lemma \ref{Lemma 4.5}, there holds
$$\xi_n=(u^*_n)_{\beta_n}=\left[(u_n)_{\beta_n}\right]^*,
\quad\eta_n=(v^*_n)_{\beta_n}=\left[(-v_n)_{\beta_n}\right]^*.$$
Hence, note that $(4.7)^b$, $Q(\xi_n,-\eta_n)=Q\left((u^*_n)_{\beta_n},-(v^*_n)_{\beta_n}\right)
=Q\left(\left[(u_n)_{\beta_n}\right]^*,-\left[(-v_n)_{\beta_n}\right]^*\right)$ can be formulated as
\begin{equation*}
\left.
\begin{array}{ll}
Q(\xi_n,-\eta_n)&=Q\left(\left[(u_n)_{\beta_n}\right]^*,-\left[(-v_n)_{\beta_n}\right]^*\right)
\\[0.5cm]
&=\displaystyle\dfrac{a_2}{m_1}\int_{\mathbb{R}^N}|\nabla \xi_n|^2dx+\dfrac{a_1}{2m_2}\int_{\mathbb{R}^N}|\nabla \eta_n|^2dx
-\dfrac{N}{2}a_1a_2\int_{\mathbb{R}^N}  \eta_n {\xi_n}^2dx
 \\[0.5cm]
&=0.
\end{array}
\right.\eqno(4.8)^*
\end{equation*}
Note that \eqref{1.10}, $4< N<6$, by $(4)$ and $(5)$ of Lemma \ref{Lemma 4.5}, as well as \eqref{4.3}, in view of $Q(u_n,v_n)=0$ and Lemma \ref{Lemma 4.7} for $\beta_{n}=1$, we obtain
\begin{equation}\label{4.9}
\left.
\begin{array}{ll}
\displaystyle\mathcal{S}(\xi_n,-\eta_n)\leq \mathcal{S}\left((u_n)_{\beta_n},-(-v_n)_{\beta_n}\right)\leq \mathcal{S}(u_n,v_n).
\end{array}
\right.
\end{equation}
Thus \eqref{4.8} and \eqref{4.9} yield that
\begin{equation}\label{4.10}
\left.
\begin{array}{ll}
(\xi_n,-\eta_n)\in \mathcal{M} \quad and\quad\mathcal{S}(\xi_n,-\eta_n)\leq \mathcal{S}(u_n,v_n).
\end{array}
\right.
\end{equation}
This implies that $\left\{(\xi_n,-\eta_n), n\in\mathbb{N}\right\}$ itself is a minimizing sequence for \eqref{1.9}.\\
\indent Now for the minimizing sequence $\{(\xi_n,-\eta_n), n\in\mathbb{N}\}$ of \eqref{1.9}, by Lemma \ref{Lemma 4.6} and \eqref{4.6}, one gets that $\|\xi_n\|_{H^1(\mathbb{R}^N)}$ and $\|\eta_n\|_{H^1(\mathbb{R}^N)}$ are both bounded for all $n\in\mathbb{N}$. Recall that $(\xi_n,\eta_n)$ are sequences of spherically symmetric non-increasing functions, there exists a subsequence $\left\{\xi_{n_k}, k\in\mathbb{N}\right\}\subset\{\xi_{n}, n\in\mathbb{N}\}$ such that as $k\rightarrow+\infty$,
\begin{equation}\label{4.11}
\left\{
\begin{array}{ll}
\displaystyle\xi_{n_k}\rightharpoonup\xi_{\infty}\quad weakly \quad in \quad H^1(\mathbb{R}^N),
\\[0.3cm]
\xi_{n_k}\rightarrow\xi_{\infty}\quad a.e. \quad in \quad \mathbb{R}^N.
\end{array}
\right.
\end{equation}
On the other hand, for $\left\{\eta_{n_k}, k\in\mathbb{N}\right\}\subset\{\eta_{n}, n\in\mathbb{N}\}$, there also exists a subsequence $\left\{\eta_{n_{k_m}}, m\in\mathbb{N}\right\}\subset\{\eta_{n_k}, k\in\mathbb{N}\}$ such that as $m\rightarrow+\infty$
\begin{equation}\label{4.12}
\left\{
\begin{array}{ll}
\displaystyle\eta_{n_{k_m}}\rightharpoonup\eta_{\infty}\quad weakly \quad in \quad H^1(\mathbb{R}^N),
\\[0.3cm]
\eta_{n_{k_m}}\rightarrow\eta_{\infty}\quad a.e. \quad in \quad \mathbb{R}^N.
\end{array}
\right.
\end{equation}
\eqref{4.11} also yields that as $m\rightarrow+\infty$,
\begin{equation}\label{4.13}
\left\{
\begin{array}{ll}
\displaystyle\xi_{n_{k_m}}\rightharpoonup\xi_{\infty}\quad weakly \quad in \quad H^1(\mathbb{R}^N),\\\\
\xi_{n_{k_m}}\rightarrow\xi_{\infty}\quad a.e. \quad in \quad \mathbb{R}^N.
\end{array}
\right.
\end{equation}
Thus, we can extract a subsequence $\left\{\left(\xi_{n_{k_m}},\eta_{n_{k_m}}\right):~ m\in\mathbb{N}\right\}$ from
$\left\{(\xi_{n},\eta_{n}): n\in\mathbb{N}\right\}$ such that \eqref{4.12} and \eqref{4.13} hold. Without any confusion, we still label $\left\{\left(\xi_{n_{k_m}},\eta_{n_{k_m}}\right):~m\in\mathbb{N}\right\}$ with $\left\{(\xi_{n},\eta_{n}):n\in\mathbb{N}\right\}$.\\
Note that Lemma 4.8 that for $1<p<\dfrac{2N}{N-2}$, the embedding $H^1_r(\mathbb{R}^N)\hookrightarrow L^p_r(\mathbb{R}^N)$ is compact, it follows from \eqref{4.12} and \eqref{4.13} that
\begin{equation}\label{4.14}
\left\{
\begin{array}{ll}
\xi_{n}\rightarrow\xi_{\infty}\quad strongly \quad in \quad L^{2q}(\mathbb{R}^N),\\\\
\eta_{n}\rightarrow\eta_{\infty}\quad strongly \quad in \quad L^{p}(\mathbb{R}^N),
\end{array}
\right.
\end{equation}
where
$$\dfrac{1}{p}+\frac{1}{q}=1,~ 4<N<6,~ 2+\dfrac{4}{N}<p<\dfrac{2N}{N-2},~ \dfrac{2N}{N+2}<q<\frac{2N+4}{N+4}.\eqno(4.14)^*$$
We then claim:\quad  {\bf $\dfrac{a_2}{m_1}||\nabla\xi_n||^2_{L^2(\mathbb{R}^N)}+\dfrac{a_1}{2m_2}||\nabla\eta_n||^2_{L^2(\mathbb{R}^N)}$ is bounded away from $0$.}
\\
\indent Indeed, recalling $(4.7)^b$ and $(4.8)^*$, $Q(\xi_n,-\eta_n)=0$ implies that
\begin{equation}\label{4.15}
\left.
\begin{array}{ll}
\displaystyle\dfrac{a_2}{m_1}\int_{\mathbb{R}^N}|\nabla \xi_n|^2dx+\dfrac{a_1}{2m_2}\int_{\mathbb{R}^N}|\nabla \eta_n|^2dx-\dfrac{N}{2}a_1a_2\int_{\mathbb{R}^N}\eta_n\xi_n^2dx=0.
\end{array}
\right.
\end{equation}
Note that assumption \eqref{1.10}, and $4<N<6$, by Young's inequality and the Gagliardo-Nirenberg inequality, for $\dfrac{1}{p}+\dfrac{1}{q}=1$, we have
\begin{equation}\label{4.16}
\left.
\begin{array}{ll}
\displaystyle\dfrac{N}{2}a_1a_2\int_{\mathbb{R}^N}\eta_n\xi_n^2dx&
\displaystyle\leq C\left(\int_{\mathbb{R}^N}\dfrac{|\eta_n|^p}{p}dx
+\int_{\mathbb{R}^N}\dfrac{|\xi_n|^{2q}}{q}dx\right)
\\[0.5cm]
&\displaystyle\leq C\left(\left\|\eta_n\right\|^p_{L^p(\mathbb{R}^N)}+\left\|\xi_n\right\|^{2q}_{L^{2q}(\mathbb{R}^N)}\right)
\\[0.5cm]
&\displaystyle\leq C\left\|\xi_n\right\|^{2q-\frac{N}{2}(2q-2)}_{L^{2}(\mathbb{R}^N)}
\left\|\nabla\xi_n\right\|^{\frac{N}{2}(2q-2)}_{L^{2}(\mathbb{R}^N)}
\\[0.5cm]
&\displaystyle~~+C\left\|\eta_n\right\|^{p-\frac{N}{2}(p-2)}_{L^{2}(\mathbb{R}^N)}
\left\|\nabla\eta_n\right\|^{\frac{N}{2}(p-2)}_{L^{2}(\mathbb{R}^N)}.\\
\end{array}
\right.
\end{equation}
Recalling $(4.14)^*$, we obtain
\begin{equation}\label{4.17}
\left\{
\begin{array}{ll}
\dfrac{N}{4}\left(p-2\right)>1,\quad \dfrac{4N}{N+2}<2q<\dfrac{4N+8}{N+4},
\\[0.5cm]
\dfrac{N}{4}\left(2q-2\right)>\dfrac{N}{4}\left(\dfrac{4N}{N+2}-2\right)
=\dfrac{N}{2}\dfrac{N-2}{N+2}>1.
\end{array}
\right.
\end{equation}
Let $\theta=max\left\{\dfrac{N}{4}(p-2),\dfrac{N}{4}(2q-2)\right\}>1$, recalling that $||\xi_n||_{L^2(\mathbb{R}^N)}\leq C$, $||\eta_n||_{L^2(\mathbb{R}^N)}\leq C$, \eqref{4.16} and \eqref{4.17} yield that
\begin{equation*}
\left.
\begin{array}{ll}
\displaystyle\dfrac{a_2}{m_1}\int_{\mathbb{R}^N}|\nabla \xi_n|^2dx+\dfrac{a_1}{2m_2}\int_{\mathbb{R}^N}|\nabla \eta_n|^2dx\leq c\left(\dfrac{a_2}{m_1}\int_{\mathbb{R}^N}|\nabla \xi_n|^2dx+\dfrac{a_1}{2m_2}\int_{\mathbb{R}^N}|\nabla \eta_n|^2dx\right)^{\theta},
\end{array}
\right.
\end{equation*}
which implies that
\begin{equation*}
\left.
\begin{array}{ll}
\displaystyle\left(\dfrac{a_2}{m_1}\|\nabla\xi_n\|^2
_{L^2(\mathbb{R}^N)}+\dfrac{a_1}{2m_2}\|\nabla\eta_n\|^2_{L^2(\mathbb{R}^N)}\right)^{\theta-1}\geq C>0,
\end{array}
\right.
\end{equation*}
that is,
\begin{equation}\label{4.18}
\left.
\begin{array}{ll}
\displaystyle\dfrac{a_2}{m_1}\|\nabla\xi_n\|^2_{L^2(\mathbb{R}^N)}
+\dfrac{a_1}{2m_2}\|\nabla\eta_n\|^2_{L^2(\mathbb{R}^N)}\geq C>0.
\end{array}
\right.
\end{equation}
This gives that $\dfrac{a_2}{m_1}||\nabla\xi_n||^2_{L^2(\mathbb{R}^N)}+\dfrac{a_1}{2m_2}||\nabla\eta_n||^2_{L^2(\mathbb{R}^N)}$ is bounded away from $0$. Thus, we claim:
\begin{equation*}
\left.
\begin{array}{ll}
(\xi_{\infty},\eta_{\infty})\neq(0,0).
\end{array}
\right.\eqno(4.18)^*
\end{equation*}
\indent Indeed, by contradiction, if $(\xi_{\infty},\eta_{\infty})\equiv(0,0)$, then for $$\frac{1}{p}+\frac{1}{q}=1, \quad 2+\frac{4}{N}<p<\frac{2N}{N-2},\quad \frac{2N}{N+2}<q<\frac{2N}{N+4}, \quad 4<N<6,$$
we obtain
\begin{equation*}
\left\{
\begin{array}{ll}
\xi_{n}\rightarrow0\quad strongly \quad in \quad L^{2q}(\mathbb{R}^N),\\\\
\eta_{n}\rightarrow0\quad strongly \quad in \quad L^{p}(\mathbb{R}^N).
\end{array}
\right.
\end{equation*}
Since by \eqref{4.15}, as $n\rightarrow+\infty$,
\begin{equation*}
\left.
\begin{array}{ll}
0<c&\displaystyle\leq\dfrac{a_2}{m_1}\int_{\mathbb{R}^N}|\nabla \xi_n|^2dx+\dfrac{a_1}{2m_2}\int_{\mathbb{R}^N}|\nabla \eta_n|^2dx
\\[0.5cm]
&\displaystyle=\dfrac{N}{2}a_1a_2\int_{\mathbb{R}^N}\eta_n\xi_n^2dx
\\[0.5cm]
&\leq \displaystyle c\int_{\mathbb{R}^N}|\xi_n|^{2q}dx+c\int_{\mathbb{R}^N}|\eta_n|^{p}dx\rightarrow0,
\end{array}
\right.
\end{equation*}
which implies$$\frac{a_2}{m_1}\int_{\mathbb{R}^N}|\nabla \xi_n|^2dx+\frac{a_1}{2m_2}\int_{\mathbb{R}^N}|\nabla \eta_n|^2dx\rightarrow0\,\quad as\quad n\rightarrow+\infty.$$
This contradicts to \eqref{4.18}, and hence $(4.18)^*$ holds true: $(\xi_{\infty},\eta_{\infty})\neq(0,0)$. \\
\indent Next, we let
\begin{equation}\label{4.19}
\left.
\begin{array}{ll}
\xi=(\xi_{\infty})_{\beta},\quad \eta=(\eta_{\infty})_{\beta},
\end{array}
\right.
\end{equation}
with $\beta>0$ uniquely determined by the condition
\begin{equation*}
\left.
\begin{array}{ll}
Q(\xi,-\eta)=Q((\xi_{\infty})_{\beta},-(\eta_{\infty})_{\beta})=0.
\end{array}
\right.
\end{equation*}
Here, as (4.15), $\mathcal{Q}(\xi,-\eta)=\mathcal{Q}((\xi_{\infty})_{\beta}, -(\eta_{\infty})_{\beta})$ can be expressed by
\begin{equation*}
\left.
\begin{array}{ll}
\displaystyle \mathcal{Q}(\xi,-\eta)=\dfrac{a_2}{m_1}\int_{\mathbb{R}^N}|\nabla \xi |^2dx+\dfrac{a_1}{2m_2}\int_{\mathbb{R}^N}|\nabla \eta |^2dx-\dfrac{N}{2}a_1a_2\int_{\mathbb{R}^N}\eta \xi ^2dx.
\end{array}
\right.\eqno(4.19)^a
\end{equation*}
Similarly, $\mathcal{S}(\xi,-\eta)=\mathcal{S}((\xi_{\infty})_{\beta}, -(\eta_{\infty})_{\beta})$  can be formulated as
\begin{equation*}
\left.
\begin{array}{ll}
\mathcal{S}(\xi,-\eta)=&\displaystyle\dfrac{a_2}{2m_1}\int_{\mathbb{R}^N}|\nabla \xi|^2dx+\dfrac{a_1}{4m_2}\int_{\mathbb{R}^N}|\nabla \eta|^2dx
\\[0.5cm]
&\displaystyle+a_2\omega_1\int_{\mathbb{R}^N}|\xi|^2dx+\dfrac{a_1}{2}
\omega_2\int_{\mathbb{R}^N}|\eta|^2dx-
a_1a_2\int_{\mathbb{R}^N}\eta \xi^2dx.
\end{array}
\right.\eqno(4.19)^b
\end{equation*}
Therefore we have
\begin{equation}\label{4.20}
\left\{
\begin{array}{ll}
(\xi_n)_{\beta}\rightarrow\xi\quad strongly\quad in\quad L^{2q}(\mathbb{R}^N),\\\\
(\eta_n)_{\beta}\rightarrow\eta\quad strongly\quad in\quad L^{p}(\mathbb{R}^N),\\\\
(\xi_n)_{\beta}\rightarrow\xi,\quad (\eta_n)_{\beta}\rightarrow\eta \quad weakly\quad in\quad H^{1}(\mathbb{R}^N),\\\\
(\xi_n)_{\beta}\rightarrow\xi,\quad (\eta_n)_{\beta}\rightarrow\eta \quad a.e.\quad in\quad  \mathbb{R}^N,
\end{array}
\right.
\end{equation}
where $p,q$ are decided by \eqref{4.14} and $(4.14)^*$.\\
\indent Since by \eqref{4.8}, $\mathcal{Q}(\xi_n,-\eta_n)=0$, Lemma \ref{Lemma 4.7} then yields that
\begin{equation}\label{4.21}
\left.
\begin{array}{ll}
\mathcal{S}\left((\xi_n)_{\beta},-(\eta_n)_{\beta}\right)\leq\mathcal{S}
\left(\xi_n,-\eta_n\right)
\end{array}
\right.
\end{equation}
\eqref{4.20} and \eqref{4.21} then imply that
\begin{equation}\label{4.22}
\left.
\begin{array}{ll}
\mathcal{S}\left(\xi,-\eta\right)&\leq\liminf\limits_{n\rightarrow+\infty}
\mathcal{S}\left[(\xi_n)_{\beta},-(\eta_n)_{\beta}\right]
\\[0.5cm]
&\leq\lim\limits_{n\rightarrow+\infty}\mathcal{S}\left(\xi_n,-\eta_n\right)
=\inf\limits_{(u,v)\in\mathcal{M}}\mathcal{S}(u,v).
\end{array}
\right.
\end{equation}
Note that $(\xi,\eta)\neq(0,0)$ and $\mathcal{Q}(\xi,-\eta)=0$, there holds $(\xi,-\eta)\in\mathcal{M}$. Therefore, \eqref{4.22} yields that $(\xi,-\eta)$ solves the minimization problem:
\begin{equation}\label{4.23}
\left.
\begin{array}{ll}
\mathcal{S}\big(\xi,-\eta\big)=\min\limits_{(u,v)\in\mathcal{M}}\mathcal{S}(u,v).
\end{array}
\right.
\end{equation}
This completes the proof of Proposition \ref{Proposition 4.3}.
\end{proof}
So far, the proof of $(1)$ in Theorem \ref{Theorem 4.1} is finished.
\end{proof}
\subsection{Existence of ground state solution of \eqref{4.1}
 \\[0.3cm]  \indent ----proof of  (2) in Theorem  \ref{Theorem 4.1}}
\indent~~In this subsection, we prove $(2)$ of Theorem \ref{Theorem 4.1}.
 \begin{proof}Since $(\xi,-\eta)$ is a solution of the minimization problem \eqref{4.23}, there exists a Lagrange multiplier $\Lambda$ such that
\begin{equation}\label{4.24}
\left.
\begin{array}{ll}
\delta_{\xi}\left[\mathcal{S}(\xi,-\eta)+\Lambda \mathcal{Q}(\xi,-\eta)\right]=0,\quad\delta_{-\eta}\left[\mathcal{S}(\xi,-\eta)+\Lambda \mathcal{Q}(\xi,-\eta)\right]=0,
\end{array}
\right.
\end{equation}
where $\delta_{u}T$ denotes the variation of $T(u,v)$ about $u$. Note that the formula
$$\delta_{u}T(u,v)=\dfrac{\partial}{\partial\zeta}T\left(u+\zeta\delta u,v\right)\Big|_{\zeta=0},$$
one has
\begin{equation}\label{4.25}
\left\{
\begin{array}{ll}
\delta_{\xi}\left[\mathcal{S}(\xi,-\eta)+\Lambda \mathcal{Q}(\xi,-\eta)\right]=\left\langle\mathcal{A}(\xi,-\eta),\delta\xi\right\rangle,
\\[0.5cm]
\delta_{-\eta}\left[\mathcal{S}(\xi,-\eta)+\Lambda \mathcal{Q}(\xi,-\eta)\right]=\left\langle\mathcal{B}(\xi,-\eta),\delta(-\eta)\right\rangle,
\end{array}
\right.
\end{equation}
where $\delta u$ denotes the variation of $u$, $\displaystyle\langle f,g\rangle=\int_{\mathbb{R}^N}fgdx$,
\begin{equation}\label{4.26}
\left\{
\begin{array}{ll}
\mathcal{A}(\xi,-\eta)=2(1+2\Lambda)\dfrac{a_2}{2m_1}(-\Delta\xi)+2a_2\omega_1\xi
-2\left(1+\dfrac{N}{2}\Lambda\right)a_1a_2\xi\eta,
\\[0.5cm]
\mathcal{B}(\xi,-\eta)=2(1+2\Lambda)\dfrac{a_1}{4m_2} \Delta\eta
-a_1\omega_2\eta+2\left(1+\dfrac{N}{2}\Lambda\right)a_1a_2\xi^2.
\end{array}
\right.
\end{equation}
Combining \eqref{4.24} with \eqref{4.25} and \eqref{4.26} yields
\begin{equation}\label{4.27}
\left.
\begin{array}{ll}
\displaystyle(1+2\Lambda)\dfrac{a_2}{2m_1}\int_{\mathbb{R}^N}|\nabla\xi|^2dx
+a_2\omega_1\int_{\mathbb{R}^N}\xi^2dx-\left(1+\dfrac{N}{2}\Lambda\right)
a_1a_2\int_{\mathbb{R}^N}\xi^2\eta dx=0,
\end{array}
\right.
\end{equation}
\\
\begin{equation}\label{4.28}
\left.
\begin{array}{ll}
\displaystyle(1+2\Lambda)\dfrac{a_1}{4m_2}\int_{\mathbb{R}^N}|\nabla\eta|^2dx
+\dfrac{a_1}{2}\omega_2\int_{\mathbb{R}^N}\eta^2dx-\dfrac{1}{2}
\left(1+\dfrac{N}{2}\Lambda\right)a_1a_2\int_{\mathbb{R}^N}\xi^2\eta dx=0.
\end{array}
\right.
\end{equation}
\\
On the other hand $Q(\xi,-\eta)=0$ gives
\begin{equation}\label{4.29}
\left.
\begin{array}{ll}
\displaystyle 2\left(\dfrac{a_2}{2m_1}\int_{\mathbb{R}^N}|\nabla \xi|^2dx+\dfrac{a_1}{4m_2}\int_{\mathbb{R}^N}|\nabla \eta|^2dx\right)-\dfrac{N}{2}a_1a_2\int_{\mathbb{R}^N}\eta \xi^2dx=0,
\end{array}
\right.
\end{equation}
which is equivalent to
\begin{equation}\label{4.30}
\left.
\begin{array}{ll}
\displaystyle \dfrac{6}{N}\left(\dfrac{a_2}{2m_1}\int_{\mathbb{R}^N}|\nabla \xi|^2dx+\dfrac{a_1}{4m_2}\int_{\mathbb{R}^N}|\nabla \eta|^2dx\right)-\dfrac{3}{2}a_1a_2\int_{\mathbb{R}^N}\eta \xi^2dx=0.
\end{array}
\right.
\end{equation}
Combining \eqref{4.27}, \eqref{4.28}, \eqref{4.29} and \eqref{4.30} together yields
\begin{equation}\label{4.31}
\left.
\begin{array}{ll}
&\displaystyle \left(1-\dfrac{6}{N}\right)\left(\dfrac{a_2}{2m_1}\int_{\mathbb{R}^N}|\nabla \xi|^2dx+\dfrac{a_1}{4m_2}\int_{\mathbb{R}^N}|\nabla \eta|^2dx\right)
\\[0.5cm]
&\qquad+ \displaystyle a_2\omega_1\int_{\mathbb{R}^N}\xi^2dx+\dfrac{a_1}{2}\omega_2\int_{\mathbb{R}^N}\eta^2dx
-
\dfrac{N}{4}\Lambda a_1a_2\int_{\mathbb{R}^N}\eta\xi^2dx=0.
\end{array}
\right.
\end{equation}
Let\begin{equation}\label{4.32}
 \xi^{\tau}(x)=\frac{1}{\tau^{2}}\xi\left(\frac{x}{\tau}\right),\quad  \eta^{\tau}(x)=\frac{1}{\tau^{2}}\eta\left(\frac{x}{\tau}\right), \quad\tau>0,
 \end{equation}
 then
\begin{equation}\label{4.33}
\left.
\begin{array}{ll}
\displaystyle &Q(\xi^{\tau},-\eta^{\tau})=\tau^{N-6}\dfrac{a_2}{m_1}\int_{\mathbb{R}^N}|\nabla \xi|^2dx+\tau^{N-6}\dfrac{a_1}{2m_2}\int_{\mathbb{R}^N}|\nabla \eta|^2dx-\frac{N}{2}\tau^{N-6}a_1a_2\int_{\mathbb{R}^N}\eta\xi^2dx
\\[0.5cm]
&\qquad \tau^{N-6}\mathcal{Q}(\xi,-\eta).
\end{array}
\right.
\end{equation}
Since $(\xi,-\eta)\in\mathcal{M}$, \eqref{1.7} and \eqref{4.33} imply that
$$\forall \tau>0,\quad \left(\xi^{\tau},-\eta^{\tau}\right)\in\mathcal{M}.\eqno(4.33a)$$
By Lemma \ref{Lemma 4.7}, it follows that the function $\tau\rightarrow\mathcal{S}\left(\xi^{\tau},-\eta^{\tau}\right)$ attains a minimum at $\tau=1$, which yields that
\begin{equation}\label{4.34}
\left.
\begin{array}{ll}
\dfrac{d}{d\tau}\mathcal{S}\left(\xi^{\tau},-\eta^{\tau}\right)\Big|_{\tau=1}=0.
\end{array}
\right.
\end{equation}
From \eqref{4.3} and $\left(\xi^{\tau},-\eta^{\tau}\right)\in\mathcal{M}$, $\mathcal{S}\left(\xi^{\tau},-\eta^{\tau}\right)$ has the following expression
\begin{equation}\label{4.35}
\left.
\begin{array}{ll}
\displaystyle \mathcal{S}\left(\xi^{\tau},-\eta^{\tau}\right)&\displaystyle
=\left(1-\dfrac{4}{N}\right)\tau^{N-6}\left(\dfrac{a_2}{2m_1}\int_{\mathbb{R}^N}|\nabla \xi|^2dx+\dfrac{a_1}{4m_2}\int_{\mathbb{R}^N}|\nabla \eta|^2dx\right)
\\[0.5cm]
&\displaystyle~~+\tau^{N-4}\left(a_2\omega_1\int_{\mathbb{R}^N}\xi^2dx
+\dfrac{a_1}{2}\omega_2\int_{\mathbb{R}^N}\eta^2dx\right).
\end{array}
\right.
\end{equation}
\eqref{4.34} and \eqref{4.35} then yield
\begin{equation*}
\left.
\begin{array}{ll}
&\displaystyle\left(1-\frac{4}{N}\right)(N-6)\left(\dfrac{a_2}{2m_1}
\int_{\mathbb{R}^N}|\nabla \xi|^2dx+\dfrac{a_1}{4m_2}\int_{\mathbb{R}^N}|\nabla \eta|^2dx\right)
\\[0.5cm]
&~~\displaystyle+(N-4)\left(a_2\omega_1\int_{\mathbb{R}^N}\xi^2dx
+\dfrac{a_1}{2}\omega_2\int_{\mathbb{R}^N}\eta^2dx\right)=0,
\end{array}
\right.
\end{equation*}
that is,
\begin{equation}\label{4.36}
\left.
\begin{array}{ll}
&\displaystyle\dfrac{N-6}{N}\left(\dfrac{a_2}{2m_1}\int_{\mathbb{R}^N}|\nabla \xi|^2dx+\dfrac{a_1}{4m_2}\int_{\mathbb{R}^N}|\nabla \eta|^2dx\right)
\\[0.5cm]
&\qquad\displaystyle+a_2\omega_1\int_{\mathbb{R}^N}\xi^2dx+\dfrac{a_1}{2}\omega_2\int_{\mathbb{R}^N}\eta^2dx=0.
\end{array}
\right.
\end{equation}
Recalling \eqref{4.31}, \eqref{4.36} yields
\begin{equation}\label{4.37}
\left.
\begin{array}{ll}
\displaystyle\dfrac{N}{4}\Lambda a_1a_2\int_{\mathbb{R}^N}\eta\xi^2dx=0.
\end{array}
\right.
\end{equation}
On the other hand, noting that $\mathcal{Q}(\xi,-\eta)=0$ and $(\xi,\eta)\neq(0,0)$, \eqref{4.37} and \eqref{1.7} imply
\begin{equation}\label{4.38}
\left.
\begin{array}{ll}
\displaystyle\frac{1}{2}\Lambda\left(\dfrac{a_2}{m_1}\int_{\mathbb{R}^N}|\nabla \xi|^2dx+\dfrac{a_1}{2m_2}\int_{\mathbb{R}^N}|\nabla \eta|^2dx\right)=0,
\end{array}
\right.
\end{equation}
which claims that $\Lambda=0$.  Hence,  \eqref{4.25} and \eqref{4.26} then give
\begin{equation*}
\left\{
\begin{array}{ll}
-\dfrac{a_2}{2m_1}\Delta\xi+a_2\omega_1\xi=a_1a_2\xi\eta,\\\\
-\dfrac{a_1}{2m_2}\Delta\eta+a_1\omega_2\eta=a_1a_2\xi^2,
\end{array}
\right.
\end{equation*}
which together with $a_1>0, a_2>0$ is equivalent to
\begin{equation*}
\left\{
\begin{array}{ll}
-\dfrac{1}{2m_1}\Delta\xi+ \omega_1\xi=-a_1 \xi(-\eta),\\\\
-\dfrac{1}{2m_2}\Delta(-\eta)+ \omega_2(-\eta)=- a_2\xi^2.
\end{array}
\right.\eqno(4.38)^*
\end{equation*}
This implies that  $(\xi,-\eta)$ is a solution of \eqref{4.1}, as \eqref{4.1} is the Euler-Lagrange equations of the functional $S$ (see for \eqref{1.6}). Recalling Lemma \ref{Lemma 4.8}, there holds $(\xi,-\eta)\in\mathcal{M}$. Therefore, $(\xi,-\eta)$ is then a ground state solution of \eqref{4.1}.\\
\indent This completes the proof of $(2)$ in Theorem \ref{Theorem 4.1}.
\end{proof}
\subsection{Exponential decay of the ground state solution of \eqref{4.1}
  \\[0.3cm]  \indent ----proof of  (3) in Theorem  \ref{Theorem 4.1}}
\indent~~ In this subsection, we will show the exponential decay of $(\xi,-\eta)$, for $(\xi,-\eta)$ being a spherically symmetric solution of \eqref{4.1} at infinity under the conditions in Lemma \ref{Lemma 4.2}. We first claim:
\begin{prop}\label{Proposition 4.9} Let Lemma \ref{Lemma 4.2} hold true, and let $(\xi,-\eta)$ be the solution of \eqref{4.1} obtained in $(1)$ and $(2)$ of Theorem \ref{Theorem 4.1}. Then $(\xi,-\eta)$ satisfies the following exponential decay estimates:
\begin{equation}\label{4.39}
\left.
\begin{array}{ll}
\displaystyle \left|D^{\alpha}\xi(x)\right|\leq ce^{-\delta|x|},\quad\left|D^{\beta}(-\eta(x))\right|\leq ce^{-\kappa|x|},\quad x\in \mathbb{R}^N
\end{array}
\right.
\end{equation}
for some $c>0,~\delta>0,~\kappa>0$,~$|\alpha|\leq2$ and $|\beta|\leq2$.
\end{prop}
Before showing Proposition \ref{Proposition 4.9}, we first investigate the regularity of the ground state solution $(\xi,-\eta)$ of the problem \eqref{4.1} by using $(S-4)^*$ and $(S-5)^*$ in Lemma \ref{Lemma 4.2}. We then claim:
\begin{lem}\label{Lemma 4.10}Assuming that $\omega_2=2\omega_1, ~4<N<6$, and \eqref{1.10} holds true, let $(\xi,-\eta)$ be the solution of \eqref{4.1} obtained in $(1)$ and $(2)$ of Theorem \ref{Theorem 4.1}. Then $(\xi,-\eta)\in C^2(\mathbb{R}^N)\times C^2(\mathbb{R}^N)$.
\end{lem}
\begin{proof}
Since $(\xi,\eta)$ is a spherically symmetric solution of \eqref{4.1}, $(\xi,-\eta)$ satisfies for $\omega_2=2\omega_1$,
\begin{equation}\label{4.40}
\left\{
\begin{array}{ll}
\displaystyle-\dfrac{1}{2m_1}\Delta\xi(x)=-\omega_1\xi(x)-a_1\xi(x)(-\eta(x))=g_1(\xi,-\eta),
\\[0.3cm]
\displaystyle-\dfrac{1}{2m_2}\Delta(-\eta(x))=-\omega_2(-\eta(x))-a_2\xi(x)^2=g_2(\xi,-\eta).
\end{array}
\right.
\end{equation}
Note that
\begin{equation*}\label{4.40a}
\left\{
\begin{array}{ll}
~~\displaystyle g_1(\xi,-\eta)=\left(-\omega_1-a_1(-\eta(x))\right)
\xi(x)\triangleq\dfrac{g_1(\xi,-\eta)}{\xi(x)}\cdot\xi(x),
\\[0.5cm]
~~\displaystyle g_2(\xi,-\eta)=\left(-\omega_2-a_2\dfrac{\xi^2(x)}{-\eta(x)}\right)(-\eta(x))
\triangleq\dfrac{g_2(\xi,-\eta)}{-\eta(x)}\cdot(-\eta(x)),\\
\end{array}
\right.\eqno(4.40a)
\end{equation*}
\\
by $(S-4)^*$ and $(S-5)^*$ of Lemma \ref{Lemma 4.2} with $-\eta$ playing the role of $v$  in Lemma \ref{Lemma 4.2}, one obtains for $L=\dfrac{N+2}{N-2}$,
\begin{equation}\label{4.41}
\left\{
\begin{array}{ll}
~~\displaystyle \left|\dfrac{g_1(\xi,-\eta)}{\xi}\right|\leq c+|\xi|^{L-1}=c+|\xi|^{\frac{4}{N-2}},\qquad (4.41-a)
\\[0.5cm]
\displaystyle
~~\left|\dfrac{g_2(\xi,-\eta)}{-\eta}\right|\leq c+|\eta|^{L-1}=c+|\eta|^{\frac{4}{N-2}}.\qquad (4.41-b)
\end{array}
\right.
\end{equation}
\\
Since $(\xi,-\eta)\in H^1(\mathbb{R}^N)\times H^1(\mathbb{R}^N)$, we have $(\xi,-\eta)\in L^{2^*}(\mathbb{R}^N)\times L^{2^*}(\mathbb{R}^N)$ for $2^*=\dfrac{2N}{N-2}$. Noting that
 $2^*=\dfrac{4}{N-2}\cdot\dfrac{N}{2}$, (4.41-a) and (4.41-b) yield that
 $$\left(\dfrac{g_1(\xi,-\eta)}{\xi},\frac{g_2(\xi,-\eta)}{-\eta}\right)\in L^{\frac{N}{2}}(\mathbb{R}^N)\times L^{\frac{N}{2}}(\mathbb{R}^N).$$ Hence it is easy to obtain
 $$\left(\Delta\xi(x),\Delta(-\eta(x))\right)\in\left(L^{\frac{N}{N-2}}
 (\mathbb{R}^N), L^{\frac{N}{N-2}}( \mathbb{R}^N )\right).$$
 Recalling that $(\xi,-\eta)$ is a spherically symmetric solution of \eqref{4.1}, applying Sobolev embedding theorem, we have $(\xi,-\eta)\in L_{loc}^{p_1}(\mathbb{R}^N  )\times L_{loc}^{p_2}(\mathbb{R}^N )$ for $1\leq p_1,p_2<+\infty$. By \eqref{4.40}, we also have $(\Delta\xi,\Delta(-\eta))\in L_{loc}^{p_1}(\mathbb{R}^N)\times L_{loc}^{p_2}(\mathbb{R}^N)$ for $1\leq p_1,p_2<+\infty$. In addition, a classical boot strap argument (on balls $B_R$) then shows that $(\xi,-\eta)\in L_{loc}^{\infty}(\mathbb{R}^N )\times L_{loc}^{\infty}(\mathbb{R}^N )$ (see Lemma \ref{Lemma 2.6}). Thus by the $L^p-estimate$ \cite{Gilbarg-Trudinger}[2, Theorem 9.11], one obtains that $(\xi,-\eta)\in W_{loc}^{2,p_1}(\mathbb{R}^N )\times  W_{loc}^{2,p_2}(\mathbb{R}^N )$, for $1<p_1,p_2<\infty$. Hence Lemma \ref{Lemma 2.7} (Rellich's compactness theorem) yields that
 $$(\xi,-\eta)\in C^{1,\alpha_1}(\mathbb{R}^N)\times  C^{1,\alpha_2}(\mathbb{R}^N)$$
for $\alpha_1\in(0,1)$ and~$\alpha_2\in(0,1),$ with $0<\alpha_1\leq1-\dfrac{1}{p_1},0<\alpha_2\leq1-\dfrac{1}{p_2}$.\\
Noting that $(\xi,-\eta)$ is a spherically symmetric solution of \eqref{4.1}, by \eqref{4.40} one knows that for $\omega_2=2\omega_1$, $(\xi,-\eta)$ satisfies the equations as follows:
\begin{equation*}
\left\{
\begin{array}{ll}
-\dfrac{1}{2m_1}\xi_{rr}-\dfrac{1}{2m_1}\dfrac{N-1}{r}\xi_r
 =-\omega_1\xi+a_1\xi \eta,
\\[0.5cm]
-\dfrac{1}{2m_2} (-\eta_{rr}) -\dfrac{1}{2m_2}\dfrac{N-1}{r}(- \eta_r)   =
-\omega_2 (-\eta)-a_2\xi^2,
\end{array}
\right.\eqno(4.42)^*
\end{equation*}
\\
which is equivalent to that $(\xi,-\eta)$ satisfies
\begin{equation}\label{4.42}
\left\{
\begin{array}{ll}
-\dfrac{1}{2m_1}\xi_{rr}-\dfrac{1}{2m_1}\dfrac{N-1}{r}\xi_r
 =-\omega_1\xi+a_1\xi \eta\triangleq g^*_1(\xi,\eta),
\\[0.5cm]
-\dfrac{1}{2m_2} \eta_{rr} -\dfrac{1}{2m_2}\dfrac{N-1}{r} \eta_r   =
-\omega_2 \eta+a_2\xi^2\triangleq g^*_2(\xi,\eta),
\end{array}
\right.
\end{equation}
where we denote
$$g^*_1(\xi,\eta)=-\omega_1\xi+a_1\xi \eta, \quad g^*_2(\xi,\eta)=-\omega_2 \eta+a_2\xi^2.\eqno(4.42)^a$$
Hence $\xi_{rr}(r)$ and $\eta_{rr}(r)$ are continuous expect possible at $r=0$. We then claim :{\bf  $\xi_{rr}(r)$ and $\eta_{rr}(r)$ are also continuous at $r=0$.}

\begin{lem}\label{Lemma 4.11}
$\xi_{rr}(r)$ and $\eta_{rr}(r)$ are continuous at $r=0$.
\end{lem}
\begin{proof}
Let
$$Q_1(r)=g^*_1(\xi(r),\eta(r))=-\omega_1\xi(r)+a_1\xi(r)\eta(r),$$
\\
$$Q_2(r)=g^*_2(\xi(r),\eta(r))=-\omega_2\eta(r)+a_2\xi^2(r),$$
where $Q_1(r)$ and $Q_2(r)$ are continuous on $[0,\infty)$.~~Rewriting \eqref{4.42} as
\begin{equation}\label{4.43}
\left\{\quad
\begin{array}{ll}
\displaystyle -\dfrac{1}{2m_1}\dfrac{d}{dr}\left(r^{N-1}\xi_r\right)=r^{N-1}Q_1(r),\qquad (4.43a)
\\[0.5cm]
\displaystyle-\dfrac{1}{2m_2}\dfrac{d}{dr}\left(r^{N-1}\eta_r\right)=r^{N-1}Q_2(r),\qquad (4.43b)
\end{array}
\right.
\end{equation}
\\
and then integrating from $0$ to $r$ yield
\\
\begin{equation}\label{4.44}
\left\{\quad
\begin{array}{ll}
\displaystyle r^{N-1}\xi_r=-2m_1\int_0^rs^{N-1}Q_1(s)ds,
\\[0.5cm]
\displaystyle r^{N-1}\eta_r=-2m_2\int_0^rs^{N-1}Q_2(s)ds
,\end{array}
\right.
\end{equation}
\\
let $s=rt$, then \eqref{4.44} becomes
\begin{equation}\label{4.45}
\left\{\quad
\begin{array}{ll}
\displaystyle\xi_r=-2m_1r\int_0^1t^{N-1}Q_1(rt)dt,
\\[0.5cm]
\displaystyle
\eta_r=-2m_2r\int_0^1t^{N-1}Q_2(rt)dt.
\end{array}
\right.
\end{equation}
\\
Note that
\begin{equation*}
\left\{\quad
\begin{array}{ll}
\displaystyle\lim\limits_{r\rightarrow0}\int_0^1t^{N-1}Q_1(rt)dt=\dfrac{Q_1(0)}{N},
\\[0.5cm]
\displaystyle
\lim\limits_{r\rightarrow0}\int_0^1t^{N-1}Q_2(rt)dt=\dfrac{Q_2(0)}{N},
\end{array}
\right.\eqno(4.45)^*
\end{equation*}
\\
we obtain that $\xi_{rr}(0)$ and $\eta_{rr}(0)$ exist such that
\\
\begin{equation}\label{4.46}
\left.
\begin{array}{ll}
\displaystyle\xi_{rr}(0)=-\dfrac{ 2m_1Q_1(0)}{N},\quad\eta_{rr}(0)=-\dfrac{2m_2Q_2(0)}{N}.
\end{array}
\right.
\end{equation}
\\
{\bf Indeed, direct calculation gives
\begin{equation*}
\left.
\begin{array}{ll}
 \displaystyle\xi_{rr}(0)&=\lim\limits_{r\rightarrow0}\dfrac{\xi_{r}(r)
-\xi_{r}(0)}{r}==\displaystyle \lim\limits_{r\rightarrow0}\dfrac{-2m_1r\int_0^1t^{N-1}Q_1(rt)dt}{r}
\\[0.5cm]
&
\displaystyle =\lim\limits_{r\rightarrow0}(-2m_1)\int_0^1t^{N-1}Q_1(rt)dt=-\dfrac{ 2m_1Q_1(0)}{N},
\end{array}
\right.
\end{equation*}
and
\begin{equation*}
\left.
\begin{array}{ll}
 \displaystyle\eta_{rr}(0)&=\lim\limits_{r\rightarrow0}\dfrac{\eta_{r}(r)
-\eta_{r}(0)}{r}==\displaystyle \lim\limits_{r\rightarrow0}\dfrac{-2m_2r\int_0^1t^{N-1}Q_2(rt)dt}{r}
\\[0.5cm]
&
\displaystyle =\lim\limits_{r\rightarrow0}(-2m_2)\int_0^1t^{N-1}Q_2(rt)dt=-\dfrac{ 2m_2Q_2(0)}{N}.
\end{array}
\right.
\end{equation*}
}
On the other hand, from \eqref{4.42}, \eqref{4.45} and $(4.45)^*$ it follows that
\\
\begin{equation}\label{4.47}
\left\{\quad
\begin{array}{ll}
\lim\limits_{r\rightarrow0}\xi_{rr}(r)&=\lim\limits_{r\rightarrow0}\left
(-\dfrac{N-1}{r}\xi_r(r)-2m_1Q_1(r)\right)
\\[0.4cm]
&=-2m_1Q_1(0)+2m_1\dfrac{N-1}{N}Q_1(0)=-\dfrac{2m_1}{N}Q_1(0),
\\[0.5cm]
\lim\limits_{r\rightarrow0}\eta_{rr}(r)&=\lim\limits_{r\rightarrow0}\left
(-\dfrac{N-1}{r}\eta_r(r)-2m_2Q_2(r)\right)
=-\dfrac{2m_2}{N}Q_2(0).
\end{array}
\right.
\end{equation}
Combining \eqref{4.46} with \eqref{4.47} implies that $\xi_{rr}$ and $\eta_{rr}$ are continuous at $r=0$.\\
\indent The proof of Lemma \ref{Lemma 4.11} is then completed.
\end{proof}
So far, we obtain that $\xi_{rr}$ and $\eta_{rr}$ are continuous for any $r\geq0$, and hence Lemma \ref{Lemma 4.10} hold.
\end{proof}
We next claim $\xi>0$ and $\eta>0$ on $\mathbb{R}^N$.\\
\indent Indeed, since $\xi$ is a decreasing function of $r$, we have $\displaystyle\dfrac{d\xi}{dr}<0$ for any $r>0$. Note that $\xi(x)\in H^1(\mathbb{R}^N)$, there holds that $\lim\limits_{r\rightarrow\infty}\xi(r)=0.$ Hence by the maximum principle, one obtains that $\xi>0$ on $\mathbb{R}^N$.
On the other hand, in view of the proof of $(1)$ and $(2)$ in Theorem \ref{Theorem 4.1}, we obtain $\eta$ is a decreasing function of $r$, and hence $\displaystyle\dfrac{d\eta}{dr}<0$ for any $r>0$. From $\eta(x)\in H^1(\mathbb{R}^N)$, it follows that $\lim\limits_{r\rightarrow\infty}\eta(r)=0.$ Therefore, by the maximum principle, we get that $\eta>0$ on $\mathbb{R}^N$.
All in all, there hold $\xi(r)>0$ and $\eta(r)>0$ on $\mathbb{R}^N$.\\
\indent We are now in the position to {\bf prove Proposition \ref{Proposition 4.9}}, which will be divided into three steps.\\
\\
step1)\quad Proof of the exponential decay of $(\xi,\eta)$;
\\[0.3cm]
step2)\quad Verification of the exponential decay of $(\xi_r,\eta_r)$;
\\[0.3cm]
step3)\quad Justification of the exponential decay of $(\xi_{rr},\eta_{rr})$, and thus of $\left(|\mathcal{D}^{\alpha}\xi(x)|,|\mathcal{D}^{\alpha}\eta(x)|\right)$ for $|\alpha|\leq2$.
\begin{proof}
\quad\\
\\
{\bf step1)\quad Proof of the exponential decay of $(\xi,\eta)$}.\\
\\
\indent By Lemma \ref{Lemma 4.10}, $(\xi,\eta)$ are of class $C^2(\mathbb{R}^N)\times C^2(\mathbb{R}^N)$, accordingly it satisfies \eqref{4.42}. Set
\begin{equation}\label{4.48}
\left.
\begin{array}{ll}
\xi^*=r^\frac{N-1}{2}\xi,\quad \eta^*=r^\frac{N-1}{2}\eta.
\end{array}
\right.
\end{equation}
Direct calculation gives
\begin{equation}\label{4.49}
\left\{\quad
\begin{array}{ll}
\xi=r^{-\frac{N-1}{2}}\xi^*,\quad \eta=r^{-\frac{N-1}{2}}\eta^*,\\\\
\xi_r=\dfrac{1-N}{2}r^{-\frac{N+1}{2}}\xi^*+r^{\frac{1-N}{2}}\xi_r^*,\\\\
\xi_{rr}=r^{\frac{1-N}{2}}\xi^*_{rr}+(1-N)r^{-\frac{N+1}{2}}\xi^*_{r}+\dfrac{N^2-1}{4}r^{\frac{-N-3}{2}}\xi^*,\\\\
\eta_r=\dfrac{1-N}{2}r^{-\frac{N+1}{2}}\eta^*+r^{\frac{1-N}{2}}\eta_r^*,\\\\
\eta_{rr}=r^{\frac{1-N}{2}}\eta^*_{rr}+(1-N)r^{-\frac{N+1}{2}}\eta^*_{r}+\dfrac{N^2-1}{4}r^{\frac{-N-3}{2}}\eta^*,\\\\
\end{array}
\right.
\end{equation}
which together with \eqref{4.42} implies that $(\xi^*,\eta^*)$ satisfying
\begin{equation}\label{4.50}
\left\{\quad
\begin{array}{ll}
\xi^*_{rr}=\left[f_1(r)+\dfrac{a}{r^2}\right]\xi^*,
\\[0.5cm]
\eta^*_{rr}=\left[f_2(r)+\dfrac{a}{r^2}\right]\eta^*,
\end{array}
\right.
\end{equation}
where
\begin{equation*}f_1(r)=\dfrac{-2m_1g^*_1(\xi(r),\eta(r))}{\xi(r)},\quad f_2(r)=\dfrac{-2m_2g^*_2(\xi(r),\eta(r))}{\eta(r)},\quad a=\dfrac{(N-1)(N-3)}{4}.
\eqno(4.50)^*
\end{equation*}
\\
Recalling the radial Lemma \ref{Lemma 2.3}, $\xi(r)\rightarrow0$, $\eta(r)\rightarrow0$ as $r\rightarrow\infty$. Noting that (4.40), (4.42) and (4.42)$^a$, there holds
$$\displaystyle g^*_1\left(\xi(r),\eta(r)\right)=g_1\left(\xi(r),-\eta(r)\right), \quad g^*_2\left(\xi(r),\eta(r)\right)=-g_2\left(\xi(r),-\eta(r)\right),\eqno(4.50)^a$$
\\
this yields for $ \xi>0$ and $\eta>0$ that
\\
$$\displaystyle\dfrac{ g^*_2\left(\xi(r),\eta(r)\right)}{\eta(r)}=-\dfrac{ g_2\left(\xi(r),-\eta(r)\right)}{\eta(r)}= \dfrac{ g_2\left(\xi(r),-\eta(r)\right)}{-\eta(r)}.\eqno(4.50)^b$$
\\
Therefore from $(S-2)$ and $(S-3)$ of Lemma \ref{Lemma 4.2} and $(4.50)^*$, it follows that for $r\geq r_0$ large enough
\\
\begin{equation}\label{4.51}
\left\{
\begin{array}{ll}
f_1(r)+\dfrac{a}{r^2}\geq m_1\omega_1,
\\[0.5cm]
f_2(r)+\dfrac{a}{r^2}\geq m_2\omega_2.
\end{array}
\right.
\end{equation}
\\
Furthermore, let
$$U={\xi^*}^2,~~ V={\eta^*}^2,\eqno(4.51)^a$$
then $(U,V)$ verifies
\begin{equation}\label{4.52}
\left\{\quad
\begin{array}{ll}
\dfrac{1}{2}U_{rr}={\xi_r^*}^2+\left[f_1(r)+\dfrac{a}{r^2}\right]U,
\\[0.5cm]
\dfrac{1}{2}V_{rr}={\eta_r^*}^2+\left[f_2(r)+\dfrac{a}{r^2}\right]V,
\end{array}
\right.
\end{equation}
where we have used the fact that
 $${\xi_r^*}^2=\frac{1}{4}\frac{U_r^2}{U}\quad and \quad{\eta_r^*}^2=\frac{1}{4}\frac{V_r^2}{V}.$$
 Thus for $r\geq r_0$ one has
 \begin{equation}\label{4.53}
\left\{\quad
\begin{array}{ll}
U\geq0, V\geq0,
\\[0.3cm]
U_{rr}\geq2m_1\omega_1U,
\\[0.3cm]
V_{rr}\geq2m_2\omega_2V.
\end{array}
\right.
\end{equation}
Now let
\begin{equation}\label{4.54}
\left\{\quad
\begin{array}{ll}
U^*=e^{-\sqrt{2m_1\omega_1}r}\left(U_r+\sqrt {2m_1\omega_1}U\right),
\\[0.5cm]
V^*=e^{-\sqrt {2m_2\omega_2}r}\left(V_r+\sqrt{2m_2\omega_2}V\right),\\
\end{array}
\right.
\end{equation}
we have by \eqref{4.53}
\begin{equation}\label{4.55}
\left\{\quad
\begin{array}{ll}
U_r^*=e^{-\sqrt {2m_1\omega_1}r}\left(U_{rr}- 2m_1\omega_1 U\right)\geq0,
\\[0.5cm]
V_r^*=e^{-\sqrt[]{2m_2\omega_2}r}\left(V_{rr}- 2m_2\omega_2 V\right)\geq0.\\
\end{array}
\right.
\end{equation}
Therefore, $U^*$ and $V^*$ are nondecreasing on $(r_0,+\infty)$.\\
\indent Furthermore, we claim:
\begin{con}\label{Conclusion 4.12}
There holds
\begin{equation}\label{4.56}
\left.
\begin{array}{ll}
U^*(r)\leq0,\quad V^*(r)\leq0\quad for \quad r\geq r_1>r_0.
\end{array}
\right.
\end{equation}
\end{con}
{\bf Proof.}   Indeed, if there exists $r_1>r_0$ such that $U^*(r_1)>0$, then $U^*(r)\geq U^*(r_1)>0$ for all $r\geq r_1$. This together with \eqref{4.54} implies that
$$U_r+\sqrt{2m_1\omega_1}U=U^*(r)e^{\sqrt{2m_1\omega_1}r}\geq U^*(r_1)e^{\sqrt{2m_1\omega_1}r},$$
whence $U_r+\sqrt{2m_1\omega_1}U$ is not integrable on $(r_1,+\infty)$. Recalling (4.49), ${\xi^*}^2$ and $\xi^*\xi_r^*$ are integrable near $\infty$ for $\xi\in H^1(\mathbb{R}^N)$, this together with (4.51)$^a$ yields that $U_r+U$ are integrable from $U={\xi^*}^2$ with $U_r=2\xi^*\xi_r^*$, which is a contradiction. Hence $U^*(r)\leq 0$ for $r\geq r_1$. Similar argument yields that $V^*(r)\leq 0$ for $r\geq r_1$.\\
\indent Conclusion \ref{Conclusion 4.12} then holds true.\hfill$\Box$\\
Thus for $r\geq r_1$, \eqref{4.56} and \eqref{4.54} leads to
\begin{equation}\label{4.57}
\left\{\quad
\begin{array}{ll}
\dfrac{d\left(e^{\sqrt{2m_1\omega_1}r}U\right)}{dr}=e^{2\sqrt{2m_1\omega_1}r}U^*\leq0,
\\[0.5cm]
\dfrac{d\left(e^{\sqrt{2m_2\omega_2}r}U\right)}{dr}=e^{2\sqrt{2m_2\omega_2}r}V^*\leq0.
\end{array}
\right.
\end{equation}
This implies that
\begin{equation*}
\left\{\quad
\begin{array}{ll}
e^{\sqrt{2m_1\omega_1}r}U(r)\leq e^{\sqrt{2m_1\omega_1}r_1}U(r_1),
\\[0.5cm]
e^{\sqrt{2m_2\omega_2}r}V(r)\leq e^{\sqrt{2m_2\omega_2}r_1}V(r_1),
\end{array}
\right.
\end{equation*}
that is,
$$U(r)\leq c_1(r_1)e^{-\sqrt{2m_1\omega_1}r}\quad and\quad V(r)\leq c_2(r_1)e^{-\sqrt{2m_2\omega_2}r},$$
where~~ $c_1(r_1)\triangleq e^{\sqrt{2m_1\omega_1}r_1}U(r_1)$  and $c_2(r_1)\triangleq e^{\sqrt{2m_2\omega_2}r_1}V(r_1) $.\\
\\
Recalling that
$$U(r)={\xi^*}^2(r)=\left(r^{\frac{N-1}{2}}\xi\right)^2=r^{N-1}\xi^2,\quad V(r)={\eta^*}^2(r)=\left(r^{\frac{N-1}{2}}\eta\right)^2=r^{N-1}\eta^2,$$
we then obtain for $r\geq r_1$ and for two positive constants $c_1^*(r_1)\geq\sqrt{c_1 (r_1)}$ and $c_2^*(r_1)\geq\sqrt{c_2 (r_1)}$:
\begin{equation}\label{4.58}
\left\{\quad
\begin{array}{ll}
|\xi(r)|\leq c_1^*(r_1) r^{-\frac{N-1}{2}}e^{-\frac{\sqrt{2m_1\omega_1}}{2}r},\qquad (4.58-1)\\\\
|\eta(r)|\leq c_2^*(r_1) r^{-\frac{N-1}{2}}e^{-\frac{\sqrt{2m_2\omega_2}}{2}r}.\qquad (4.58-2)
\end{array}
\right.
\end{equation}\\
{\bf step2)\quad Verification of the exponential decay of $(\xi_r,\eta_r)$}.\\
\\
\indent We here establish the exponential decay of $\xi_r$ and $\eta_r$.\\
\indent By \eqref{4.42}
one observes that  $\xi_r$ and $\eta_r$ satisfy
\begin{equation}\label{4.59}
\left\{
\begin{array}{ll}
\dfrac{d\left(r^{N-1}\xi_r\right)}{dr}=-2r^{N-1}m_1g^*_1(\xi,\eta),\qquad (4.59-1)
\\[0.5cm]
\dfrac{d\left(r^{N-1}\eta_r\right)}{dr}=-2r^{N-1}m_2g^*_2(\xi,\eta).\qquad (4.59-2)
\end{array}
\right.
\end{equation}
Note that (4.42), $(4.42)^a$ and (4.50)$^a$, ultilizing $(S-2)$ and $(S-3)$ of Lemma \ref{Lemma 4.2} as well as the exponential decay of $(\xi,\eta)$ imply that for $r$ large enough, say $r\geq r_0$, there exists $d_1>c_1>0$ and $d_2>c_2>0$ such that\\
\begin{equation}\label{4.60}
\left.
\begin{array}{ll}
c_1|\xi|\leq\left|g^*_1(\xi,\eta)\right|=\left|g_1(\xi,-\eta)\right|\leq d_1|\xi|, \quad c_2|\eta|\leq\left|g^*_2(\xi,\eta)\right|=\left|-g_2(\xi,-\eta)\right|\leq d_2|\eta|.
\end{array}
\right.
\end{equation}
\\
 Integrating $(4.59-1)$ on $(r,R)$, applying \eqref{4.58} and letting $(r,R)\rightarrow(+\infty,+\infty)$ shows that:
 $r^{N-1}\xi_r$ has a limit as $r\rightarrow+\infty$;
 this limit can only be zero from $(4.58-1)$ and \eqref{4.60}.\\
 \indent Indeed, own to $(\xi,\eta)\in H^{1}(\mathbb{R}^n)\times H^{1}(\mathbb{R}^n)$, by $(4.59-1)$ and \eqref{4.60} one has
 \begin{equation*}
\left.
\begin{array}{ll}
\lim\limits_{r\rightarrow+\infty}r^{N-1}\xi_r&=\lim\limits_{r\rightarrow+\infty}
\int_r^{+\infty}
2s^{N-1}m_1g^*_1(\xi(s),\eta(s))ds\\\\
&\leq\lim\limits_{r\rightarrow+\infty}\int_r^{+\infty}2m_1s^{N-1}d_1|\xi(s)|ds\\\\
&\leq\lim\limits_{r\rightarrow+\infty}\int_r^{+\infty}2m_1s^{N-1}d_1c_1(r_1)s^{-\frac{N-1}{2}}e^{-\sqrt{2m_1\omega_1}s}ds\\\\
&=0.
\end{array}
\right.
\end{equation*}
Furthermore, integrating $(4.59-1)$ on $(r,+\infty)$ gives
$$-r^{N-1}\xi_r=-2\int_r^{+\infty}s^{N-1}m_1g^*_1(\xi(s),\eta(s))ds,$$
that is,
\begin{equation*}
\left.
\begin{array}{ll}
r^{N-1}\xi_r&=\int_r^{+\infty}2s^{N-1}m_1g^*_1(\xi(s),\eta(s))ds\\\\
&\leq\int_r^{+\infty}2m_1d_1|\xi(s)|ds\\\\
&\leq\int_r^{+\infty}2m_1d_1c^*_1(r_1)s^{N-1}s^{-\frac{N-1}{2}}e^{-\frac{\sqrt{2m_1\omega_1}s}{2}}ds.\\\\
\end{array}
\right.
\end{equation*}
This implies that
$$|\xi_r(r)| \leq cr^{-\theta_{1}}e^{-\frac{\sqrt{2m_1\omega_1}}{2}r}~ {\mbox for} ~{\mbox some} ~ \theta_{1}>0 ,$$ and hence $\xi_r(r)$ has an exponential decay at infinity. Similarly, combining $(4.59-2)$ with \eqref{4.60} yields
 $$|\eta_r(r)| \leq cr^{-\theta_{2}}e^{-\frac{\sqrt{2m_2\omega_2}}{2}r} ~ {\mbox for} ~{\mbox some} ~ \theta_{2}>0 ,$$ and thus $\eta_r$ has also an exponential decay at infinity.\\
 \\
{\bf step3)\quad Justification of the exponential decay of $\xi_{rr}$ and $\eta_{rr}$, and thus of $\left(|D^{\alpha}\xi(x)|,|D^{\alpha}\eta(x)|\right)$ for $|\alpha|\leq2$}.\\
\\
\indent Applying the established results that $\xi,\xi_r,\eta$ and $\eta_r$ all have exponential decays, by \eqref{4.60} and by the equivalent form of \eqref{4.42}:
\begin{equation}\label{4.61}
\left\{\quad
\begin{array}{ll}
\xi_{rr}+\dfrac{N-1}{r}\xi_r=-2m_1g^*_1(\xi,\eta),
\\[0.5cm]
\eta_{rr}+\dfrac{N-1}{r}\eta_r=-2m_2g^*_2(\xi,\eta),
\end{array}
\right.
\end{equation}
we can easily obtain the exponential decay of $\xi_{rr}$ and $\eta_{rr}$.\\
\indent By far, we have established the exponential decay of $\xi,\eta,\xi_{r},\eta_{r}, \xi_{rr}$ and $\eta_{rr}$, which imply the exponential decays for $D^{\alpha}\xi(x)$ and $D^{\alpha}\eta(x)$ for $|\alpha|\leq2.$\\
\indent Step 1), Step 2) and Step 3) then complete the proof of Proposition \ref{Proposition 4.9}.\\
\indent Finally, combining these results in {\bf Subsection 4.1, 4.2 } and {\bf 4.3} together finishes the proof of Theorem \ref{Theorem 4.1}.
\end{proof}
\section{Instability of Standing Waves under Mass Resonance }
\indent{~~} In Section 4 we have established the existence of the ground state solution $(\xi,-\eta)$ of \eqref{1.5} (or \eqref{4.1}). On the other hand, noting that
$$(\phi(t,x),\psi(t,x))=\left(e^{i\omega_1 t}\xi(x),e^{i\omega_2 t}(-\eta(x))\right)\quad with \quad\omega_2=2\omega_1,$$
which is the standing wave solution of \eqref{1.1} with the ground state solution $(\xi(x),-\eta(x))$ for \eqref{1.5} (or \eqref{4.1}). Hence Theorem \ref{Theorem 4.1} indeed gives the existence of standing waves of \eqref{1.1} associated with the ground state for \eqref{1.5} (or \eqref{4.1}). In this section, we are concerned with the characterization of the standing waves of \eqref{1.1} with minimal action $\mathcal{S}(\xi,-\eta)$ and establish its' instability.\\
\indent We then claim:
\begin{thm}\label{Theorem 5.1}\quad Let $\omega_2=2\omega_1$, $m_2=2m_1$, $4<N<6$ and \eqref{1.10} hold true. If $(\xi,-\eta)\in\mathcal{M}$ is given in Theorem \ref{Theorem 4.1}, then for arbitrary $\varepsilon>0$, there exists $(\phi_0,\psi_0)\in H^1(\mathbb{R}^N)\times H^1(\mathbb{R}^N)$ with
$$\left\|\phi_0-\xi\right\|_{H^1(\mathbb{R}^N)}<\varepsilon,\quad \left\|\psi_0-(-\eta)\right\|_{H^1(\mathbb{R}^N)}<\varepsilon$$
such that the solution $(\phi,\psi)$ of \eqref{1.1} with this initial data $(\phi_0,\psi_0)$ has the following property:\\
For some finite time $T<+\infty$, $(\phi,\psi)$ exists on $[0,T)$,
$$(\phi,\psi)\in C \left([0,T);H^1(\mathbb{R}^N)\times H^1(\mathbb{R}^N)\right),$$
and
$$\lim\limits_{t\rightarrow T}\left(||\phi||_{H^1(\mathbb{R}^N)}+||\psi||_{H^1(\mathbb{R}^N)}\right)=+\infty.$$
\end{thm}
For the local-wellposedness of solution to the Cauchy problem \eqref{1.1}-\eqref{2.1}, using the similar argument to that proposed in Cazenave \cite{Caze1} and Ginibre-Velo \cite{GV1}, we can obtain:
\begin{prop}\label{Proposition 5.2}
For any $(\phi_0,\psi_0)\in H^1(\mathbb{R}^N)\times H^1(\mathbb{R}^N)$, there exists a unique solution $(\phi,\psi)$ of \eqref{1.1}-\eqref{2.1} defined on a maximal time interval $[0,T)$, where $T=T_{max}(\phi_0,\psi_0)$ and $$(\phi,\psi)\in C \left([0,T);H^1(\mathbb{R}^N)\times H^1(\mathbb{R}^N)\right)$$ and either $T=+\infty$ or $T<+\infty$ and
$$\lim\limits_{t\rightarrow T}\left(\left\|\phi(t,\cdot)\right\|^2_{H^1(\mathbb{R}^N)}+\left\|\psi(t,\cdot)\right\|^2_{H^1(\mathbb{R}^N)}\right)=+\infty.$$
\end{prop}
\begin{prop}\label{Proposition 5.3}
Let $\omega_2=2\omega_1$, $m_2=2m_1$, $4<N<6$ and \eqref{1.10} hold true. Then for all $t\in[0,T)$, the solution $(\phi(t),\psi(t))=(\phi(t,\cdot),\psi(t,\cdot))$ of the Cauchy problem \eqref{1.1} with the initial data $(\phi_0,\psi_0)$ admits the following two conservation laws:
\begin{equation}\label{5.1}
\left.
\begin{array}{ll}
\displaystyle a_2\int_{\mathbb{R}^N}|\phi|^2dx+a_1\int_{\mathbb{R}^N}|\psi|^2dx=a_2\int_{\mathbb{R}^N}|\phi_0|^2dx+a_1\int_{\mathbb{R}^N}|\psi_0|^2dx,
\end{array}
\right.
\end{equation}\\
\begin{equation}\label{5.2}
\left.
\begin{array}{ll}
\displaystyle Re\mathcal{S}\left(\bar{\phi}(t),\psi(t)\right)=Re\mathcal{S}(\bar{\phi}_0,\psi_0).
\end{array}
\right.
\end{equation}

Put
\begin{equation}\label{5.3}
\left.
\begin{array}{ll}
\displaystyle G(t)=\int_{\mathbb{R}^N }|x|^2\left(a_2|\phi|^2+a_1|\psi|^2\right)dx,
\end{array}
\right.
\end{equation}
then
\begin{equation}\label{5.4}
\left.
\begin{array}{ll}
\displaystyle\dfrac{d^2}{dt^2}G(t)=\dfrac{2}{m_1}Re Q(\bar{\phi},\psi),
\end{array}
\right.
\end{equation}
where $\mathcal{S}(\phi,\psi)$ and $Q(\phi,\psi)$ are defined by \eqref{1.6} and \eqref{1.7}, respectively.
\end{prop}
\begin{proof}\eqref{5.1} is just the mass conservation equality \eqref{2.2}. We then verify \eqref{5.2}. From \eqref{1.6} it follows that
\begin{equation}\label{5.5}
\left.
\begin{array}{ll}
\displaystyle Re\mathcal{S}\left(\bar{\phi}(t),\psi(t)\right)&=\dfrac{a_2}{2m_1}\int_{\mathbb{R}^N}|\nabla \phi|^2dx+\dfrac{a_1}{4m_2}\int_{\mathbb{R}^N}|\nabla \psi|^2dx
\\[0.5cm]
&\displaystyle\quad +a_2\omega_1\int_{\mathbb{R}^N}|\phi|^2dx+\dfrac{a_1}{2}\omega_2
\int_{\mathbb{R}^N}|\psi|^2dx
\\[0.5cm]
&\displaystyle \quad +
a_1a_2Re\int_{\mathbb{R}^N}\psi\overline{\phi}^2dx.
\end{array}
\right.
\end{equation}
Recalling the energy conservation \eqref{2.3} and (2.3a), for $\omega_2=2\omega_1$, \eqref{5.5} gives
\begin{equation*}
\left.
\begin{array}{ll}
\displaystyle Re\mathcal{S}\left(\bar{\phi}(t),\psi(t)\right)&=\displaystyle E(\phi(t),\psi(t))+\omega_1\left[a_2\int_{\mathbb{R}^N}|\phi|^2dx+
a_1\int_{\mathbb{R}^N}|\psi|^2dx\right]
\\[0.5cm]
&\displaystyle=E(\phi_0,\psi_0)+\omega_1\left[a_2\int_{\mathbb{R}^N}|\phi_0|^2dx
+a_1\int_{\mathbb{R}^N}|\psi_0|^2dx\right]
\\[0.5cm]
&\displaystyle=Re\mathcal{S}(\bar{\phi}_0,\psi_0).
\end{array}
\right.
\end{equation*}
This proves that \eqref{5.2} holds true. In addition, by Lemma \ref{Lemma 2.2}, \eqref{1.7} and \eqref{5.3}, one obtains that for $m_2=2m_1$
\begin{equation*}
\left.
\begin{array}{ll}
\dfrac{d^2}{dt^2}G(t)&\displaystyle=\dfrac{2a_2}{m_1^2}\int_{\mathbb{R}^N}|\nabla\phi|^2dx+\dfrac{a_1}{2m_1^2}\int_{\mathbb{R}^N}|\nabla\psi|^2dx
+\dfrac{a_1a_2N}{m_1}Re\int_{\mathbb{R}^N}\psi\overline{\phi}^2dx
\\[0.5cm]
&\displaystyle=\dfrac{2}{m_1}\left[\dfrac{a_2}{m_1}\int_{\mathbb{R}^N}|\nabla\phi|^2dx
+\dfrac{a_1}{4m_1}\int_{\mathbb{R}^N}|\nabla\psi|^2dx+
\dfrac{N}{2}a_1a_2Re\int_{\mathbb{R}^N}\psi\overline{\phi}^2dx\right]
\\[0.5cm]
&\displaystyle=\dfrac{2}{m_1}Re\mathcal{Q}(\bar{\phi},\psi).
\end{array}
\right.
\end{equation*}
This completes the proof of Proposition \ref{Proposition 5.3}.
\end{proof}
{\bf Remark $5.4^*$}. If $(u,v)$ is a pair of real-valued functions, the conclusions of Lemma 4.4 are hold if we replace $\mathcal{S}(u,v)$ and $Q(u,v)$ by $Re\mathcal{S}(\bar{u},v)$ and $Re\mathcal{Q}(\bar{u},v)$, respectively. In addition, there holds $\displaystyle\int_{\mathbb{R}^N}vu^2dx=\int_{\mathbb{R}^N}v\overline{u}^2dx$.\hfill$\Box$\\

Before showing Theorem \ref{Theorem 5.1} we first claim the following:
\begin{lem}\label{Lemma 5.4}
Let $\omega_2=2\omega_1$, $m_2=2m_1$, $4<N<6$ and \eqref{1.10} hold true. Putting
\begin{equation}\label{5.6}
\left.
\begin{array}{ll}
\phi_{\mu}=\mu^{\frac{N}{2}}\phi(\mu x),\quad \psi_{\mu}=\mu^{\frac{N}{2}}\psi(\mu x),
\end{array}
\right.
\end{equation}
there exists $(\phi,\psi)\in H^1(\mathbb{R}^N)\times H^1(\mathbb{R}^N)\setminus \{(0,0)\}$ and $\mu>0$ such that
$$Re \mathcal{Q}\left(\bar{\phi}_{\mu},\psi_{\mu}\right)=0.$$
Suppose that $\mu<1$, then
\begin{equation}\label{5.7}
\left.
\begin{array}{ll}
\displaystyle Re \mathcal{S}\left(\bar{\phi},\psi\right)-Re \mathcal{S}\left(\bar{\phi}_{\mu},\psi_{\mu}\right)\geq\dfrac{1}{2}Re \mathcal{Q}(\bar{\phi},\psi).
\end{array}
\right.
\end{equation}
\end{lem}
\begin{proof}Referring to the proof of Lemma \ref{Lemma 4.7}, there holds
\begin{equation}\label{5.8}
\left\{~
\begin{array}{ll}
Re\mathcal{Q}\left(\bar{\phi}_{\lambda},\psi_{\lambda}\right)=A\lambda^2-\dfrac{N}{2}B\lambda^{\frac{N}{2}},
\\[0.5cm]
Re\mathcal{S}\left(\bar{\phi}_{\lambda},\psi_{\lambda}\right)
=A\dfrac{\lambda^2}{2}-B\lambda^{\frac{N}{2}}+C,
\end{array}
\right.
\end{equation}
 where
\begin{equation*}
\left\{~
\begin{array}{ll}
 \displaystyle A=\dfrac{a_2}{m_1}\int_{\mathbb{R}^N}|\nabla\phi|^2dx+\dfrac{a_1}{2m_2}
 \int_{\mathbb{R}^N}|\nabla\psi|^2dx,
 \\[0.5cm]
 \displaystyle B=-a_1a_2Re\int_{\mathbb{R}^N}\psi\overline{\phi}^2dx,
\\[0.5cm]
 \displaystyle C=a_2\omega_1\int_{\mathbb{R}^N}|\phi|^2dx+\dfrac{a_1}{2}\omega_2\int_{\mathbb{R}^N}|\psi|^2dx.
\end{array}
\right.
\end{equation*}
Since $m_2=2m_1>0$, $Re \mathcal{Q}\left(\bar{\phi}_{\mu},\psi_{\mu}\right)=0$ implies that
\begin{equation}\label{5.9}
\left.
\begin{array}{ll}
\displaystyle Re\int_{\mathbb{R}^N}\psi\overline{\phi}^2dx<0\quad A\mu^2=\dfrac{N}{2}B\mu^\frac{N}{2}.
\end{array}
\right.
\end{equation}
Note that $(\phi,\psi)=(\phi_1,\psi_1)$ and $Re \mathcal{Q}\left(\bar{\phi},\psi\right)=A-\dfrac{N}{2}B$, using \eqref{5.9} yields that for $0<\mu<1$ and $4<N<6$,
\begin{equation*}
\left.
\begin{array}{ll}
&\displaystyle Re\mathcal{S}\left(\bar{\phi},\psi\right)-Re\mathcal{S}\left(\bar{\phi}_{\mu},
\psi_{\mu}\right)
\\[0.5cm]
&\qquad=\dfrac{1}{2}A-B+C-\left(\dfrac{A}{2}\mu^2-B\mu^{\frac{N}{2}}+C\right)
\\[0.5cm]
&\qquad=\dfrac{1}{2}A-B-\dfrac{A}{2}\mu^2+B\mu^{\frac{N}{2}}
\\[0.5cm]
&\qquad=\dfrac{1}{2}A-\dfrac{1}{2}\cdot\dfrac{N}{2}
B\mu^{\frac{N}{2}}-B+B\mu^{\frac{N}{2}}
\\[0.5cm]
&\qquad=\dfrac{1}{2}\left(A-\dfrac{N}{2}B\right)+\dfrac{N}{4}B-\dfrac{N}{4}B\mu^{\frac{N}{2}}
-B+B\mu^{\frac{N}{2}}
\\[0.5cm]
&\qquad=\dfrac{1}{2}\left(A-\dfrac{N}{2}B\right)+\left(\dfrac{N}{4}-1\right)B
-\left(\dfrac{N}{4}-1\right)
B\mu^{\frac{N}{2}}
\\[0.5cm]
&\qquad=\dfrac{1}{2}\left(A-\dfrac{N}{2}B\right)+\left(\dfrac{N}{4}-1\right)
B\left(1-\mu^{\frac{N}{2}}\right)
\\[0.5cm]
&\qquad\geq\dfrac{1}{2}\left(A-\dfrac{N}{2}B\right)
\geq\dfrac{1}{2}Re\mathcal{Q}\left(\bar{\phi},\psi\right).
\end{array}
\right.
\end{equation*}
\end{proof}
We further establish a stronger instability result than Theorem \ref{Theorem 5.1}.\\
\indent For $(u,v)$ is a pair of complex-valued functions, we define a manifold $M^*$ as
$$M^*=\left\{(u,v)\in H^{1}(\mathbb{R}^{N})\times H^{1}(\mathbb{R}^{N}) \setminus\{(0,0)\},~~Re\mathcal{Q}\left(\bar{u},v\right)=0\right\}, \eqno(5.10)$$
and a constrained minimizing problem
$$ {M}_0=\inf\limits_{(u,v)\in\mathcal{M^*}}Re\mathcal{S}\left(\bar{u},v\right),\eqno(5.11)
$$
where $\mathcal{S}\left(\bar{u},v\right)$ and $\mathcal{Q}\left(\bar{u},v\right)$ are defined by \eqref{1.6} and \eqref{1.7} for the complex-valued pair of functions $(u,v)$, and $\bar{u}$ denotes the complex conjugate of $u$, that is,
 $$
\left.
\begin{array}{ll}
\mathcal{S}(\bar{u},v)=&\displaystyle\dfrac{a_2}{2m_1}\int_{\mathbb{R}^N}|\nabla u|^2dx+\dfrac{a_1}{4m_2}\int_{\mathbb{R}^N}|\nabla v|^2dx\\\\
&\displaystyle+a_2\omega_1\int_{\mathbb{R}^N}|u|^2dx
+\frac{a_1}{2}\omega_2\int_{\mathbb{R}^N}|v|^2dx+
a_1a_2\int_{\mathbb{R}^N}v\bar{u}^2dx,
\end{array}
\right.\eqno(5.11)^a
$$
\\
$$
\left.
\begin{array}{ll}
\mathcal{Q}(\bar{u},v)=&\displaystyle\dfrac{a_2}{m_1}\int_{\mathbb{R}^N}|\nabla u|^2dx+\dfrac{a_1}{2m_2}\int_{\mathbb{R}^N}|\nabla v|^2dx
+\dfrac{N}{2}a_1a_2\int_{\mathbb{R}^N}v \bar{u} ^2dx,
\end{array}
\right.
\eqno(5.11)^b
$$ 
 By the similar argument to that used in dealing with the real-valued pair of functions $(u,v)$, we can obtain that there exists a pair of complex-valued functions $(\xi^*,\eta^*)$ with $\displaystyle Re\int_{\mathbb{R}^{N}}\eta^*\left(\bar{\xi^*}\right)^{2}dx<0$ such that
$$Re \mathcal{S}\left(\bar{\xi^*},\eta^*\right)=M_{0}=\min\limits_{(u,v)\in M^*}
Re \mathcal{S}\left(\bar{u},v\right)>0.\eqno(5.12)$$
In addition, $\left( \xi^* ,\eta^*\right)$ satisfies the following equations:
$$
\left\{
\begin{array}{ll}
&-\dfrac{1}{2m_1}\Delta u(x)+\omega_1u(x)=-a_1v(x)\bar{u}(x),\qquad (5.12a)\\\\
&-\dfrac{1}{2m_2}\Delta v(x)+\omega_2v(x)=-a_2u^2(x),~\qquad (5.12b)
\end{array}
\right. \eqno(5.12^*)
$$
for a pair of complex-valued functions $(u,v)$, and  $\omega_{2}=2\omega_{1}$.
{\bf Remark $5.5^*$}.
Note that for a pair of real-valued functions $(u,v)$, there hold
$$
\left\{
\begin{array}{ll}
&Re \mathcal{S}\left(\bar{u},v\right)= \mathcal{S}\left(u,v\right),\quad Re \mathcal{Q}\left(\bar{u},v\right)=\mathcal{Q}(u,v),
\\[0.5cm]
&M^*~~\mbox{is~~identical~~with~~M}.
\end{array}
\right.\eqno(5.13)
$$
 From \eqref{1.9}, (5.11) and Theorem 4.1, for $(u,v)$ is a pair of real-valued functions, it follows that
 $$M_{0}=K=\mathcal{S}\left(\xi,-\eta\right).\eqno(5.14)$$
 \quad\hfill$\Box$\\
We then claim:
\begin{prop}\label{Proposition 5.5}
Let $\omega_2=2\omega_1$, $m_2=2m_1$, $4<N<6$ and \eqref{1.10} hold true. Let also $(\phi,\psi)\in H^1(\mathbb{R}^N)\times H^1(\mathbb{R}^N)$ be a solution of the Cauchy problem
\eqref{1.1}-\eqref{2.1} on $[0,T).$ 
Put\\
$$
\mathcal{K}_1=\left\{(u,v)\in H^1(\mathbb{R}^N)\times H^1(\mathbb{R}^N), ReQ\left(\bar{u},v\right)<0,~~Re\mathcal{S}
\left(\bar{u},v\right)<M_0\right\}.\eqno(5.15)
$$\\
Then for any initial data $(\phi_0,\psi_0)\in\mathcal{K}_1$, there holds $(\phi(t),\psi(t))\in\mathcal{K}_1, \forall t\in[0,T)$. That is, $\mathcal{K}_1$ is invariant under the flow generated by the Cauchy problem
\eqref{1.1}-\eqref{2.1}. Furthermore, if
$$\left(|x|\phi_0,|x|\psi_0\right)\in L^2(\mathbb{R}^N)\times L^2(\mathbb{R}^N),$$
 then $T$ is finite and
$$\lim\limits_{t\rightarrow T}\Big(||\phi(t)||_{H^1(\mathbb{R}^N)}+||\psi(t)||_{H^1(\mathbb{R}^N)}\Big)=+\infty.$$
\begin{proof}
Let $(\phi_0,\psi_0)\in\mathcal{K}_1$, since by \eqref{5.2},
$$
Re\mathcal{S}\left(\bar{\phi}(t),\psi(t)\right)=Re\mathcal{S}
\left(\bar{\phi}_0,\psi_0\right)<M_0,\quad for\quad 0\leq t<T,\eqno(5.16)
$$
we claim:
$$
Re \mathcal{Q}\left(\bar{\phi}(t),\psi(t)\right)<0 \quad for\quad 0\leq t<T.\eqno(5.17)
$$
Otherwise, by continuity there would exist a $t^*>0$ such that
$$Re\mathcal{Q}(\bar{\phi}(t^*),\psi(t^*))=0,\quad that \quad is,\quad (\phi(t^*),\psi(t^*))\in M^* ,$$
where $ M^*$ is defined by (5.10). This is impossible for
$$Re\mathcal{S}(\bar{\phi}(t^*),\psi(t^*))<{M}_0\quad and\quad  {M}_0=\min\limits_{(u,v)\in M^*}Re\mathcal{S}(\bar{u},v).$$
So (5.17) holds true. Thus, $\mathcal{K}_1$ is invariant under the flow generated by the Cauchy problem
\eqref{1.1}-\eqref{2.1}.\\
\indent Now for fixed $t\in[0,T)$, let $\mu$ be defined by $$Re\mathcal{Q}\left((\bar{\phi(t)})_{\mu},(\psi(t))_{\mu}\right)=0.$$
Note that $(\phi_0,\psi_0)\in\mathcal{K}_1$ and $\mathcal{K}_1$ is a invariant manifold, we have
 $Re\mathcal{Q}\left(\bar{\phi}(t),\psi(t)\right)<0$, which yields from Lemma \ref{Lemma 5.4} that $\mu<1$. According to  $Re\mathcal{S}\left((\bar{\phi}(t))_{\mu},(\psi(t))_{\mu}\right)\geq M_0$ and  $Re\mathcal{S}\left(\bar{\phi}(t),\psi(t)\right)=Re\mathcal{S}
\left(\bar{\phi}_0,\psi_0\right)$, by \eqref{5.7} one has
$$
\left.
\begin{array}{ll}
\displaystyle Re\mathcal{Q}(\bar{\phi}(t),\psi(t))&\displaystyle\leq2\left(Re\mathcal{S}
\left(\bar{\phi}(t),\psi(t)\right)
-Re\mathcal{S}\left((\bar{\phi}(t))_{\mu},(\psi(t))_{\mu}\right)\right)
\\[0.5cm]
&\displaystyle \leq2Re\mathcal{S}\left(\bar{\phi}_0,\psi_0\right)-2M_0<0,
\end{array}
\right.\eqno(5.18)
$$
where the last inequality uses the fact that $(\phi_0,\psi_0)\in\mathcal{K}_1$. 
Let $\theta_0=2M_0-2Re\mathcal{S}\left(\bar{\phi}_0,\psi_0\right)>0$, and $\theta_0>0$ is a fixed positive constant. Combining \eqref{5.3}, \eqref{5.4} and (5.18) together derives  that
$$
\left.
\begin{array}{ll}
& \displaystyle \dfrac{d^2}{dt^2}\left[\int_{\mathbb{R}^N}|x|^2\left(a_2|\phi|^2+a_1|\psi|^2\right)dx\right]
\\[0.5cm]
&\displaystyle \qquad=\dfrac{2}{m_1}Re\mathcal{Q}(\bar{\phi},\psi)\leq\dfrac{-2\theta_0}{m_1}<0.
\end{array}
\right.\eqno(5.19)
$$
Using the same argument as that used in the proof of Theorem \ref{Theorem 3.1}, (5.19) implies that $T$ must be finite and $$\lim\limits_{t\rightarrow T}\Big(||\phi(t)||_{H^1(\mathbb{R}^N)}+||\psi(t)||_{H^1(\mathbb{R}^N)}\Big)=+\infty.$$
This completes the proof of Proposition \ref{Proposition 5.5}.
\end{proof}
\end{prop}
Finally, we return to the proof of Theorem \ref{Theorem 5.1}.\\
\\
\indent{\bf Proof of Theorem \ref{Theorem 5.1}.}\\
\indent Put
$$
\left.
\begin{array}{ll}
\phi_0(x)=\lambda^{\frac{N}{2}}\xi(\lambda x),\quad\psi_0(x)=\lambda^{\frac{N}{2}}(-\eta)(\lambda x) \quad for \quad\lambda>1,
\end{array}
\right.\eqno(5.20)
$$
where $(\xi,\eta)$ is a pair of real-valued functions. Note that $(\xi,-\eta)\in M$ and
$$Re Q(\bar{\xi},-\eta)=Q(\xi,-\eta)=0,\eqno(5.21)$$
by (5.12), (5.13), (5.14), (5.21) and Remark 5.5$^*$,  making the similar argument to that adopted in Lemma 4.7, the functions $(\phi_0(x),\psi_0(x))$ satisfy for any $\lambda>1$,
$$Re\mathcal{Q}(\bar{\phi}_0,\psi_0)=\mathcal{Q}( \phi_0,\psi_0)<0,~~~Re\mathcal{S}(\bar{\phi}_0,\psi_0)
=\mathcal{S}( \phi_0,\psi_0)<\mathcal{S}( \xi,-\eta)=K=M_0,\eqno(5.22)$$
where $M_0$ is defined by (5.18). Hence, (5.15) and (5.22) yield that $(\phi_0,\psi_0)\in \mathcal{K}_{1}$. On the other hand, note that $\mathcal{K}_{1}$ is a invariant flow for the Cauchy problem \eqref{1.1}-\eqref{2.1} from Proposition 5.5, there holds $(\phi(t),\psi(t))\in \mathcal{K}_{1}$, where $(\phi(t),\psi(t))$ is a solution of \eqref{1.1} with the initial data $(\phi_0(x),\psi_0(x))$.\\
\indent Recalling (3) of Theorem \ref{Theorem 4.1}, $\xi(x)$ and $\eta(x)$ have the exponential decays at infinity. By (5.20), we have 
$$
\left.
\begin{array}{ll}
(|x|\phi_0,|x|\psi_0)\in L^2(\mathbb{R}^N)\times L^2(\mathbb{R}^N).
\end{array}
\right.\eqno(5.23)
$$
Furthermore, as $\lambda\rightarrow1$,
$$\left\|\phi_0-\xi\right\|_{H^1(\mathbb{R}^N)}\quad and \quad \left\|\psi_0-(-\eta) \right\|_{H^1(\mathbb{R}^N)}$$
can be chosen arbitrarily small. Thus, using Proposition \ref{Proposition 5.5}, we obtain that the solution $(\phi,\psi)$ of  \eqref{1.1} with the initial data $(\phi_0,\psi_0)$ blows up (in $H^1(\mathbb{R}^N)\times H^1(\mathbb{R}^N)-$norm) in finite time.\\
\indent This finishes the proof of Theorem 5.1.\hfill$\Box$

\section*{Acknowledgments}
Zaihui Gan is partially supported by the National Science Foundation of China under grant (No. 11571254) and the Natural Science Foundation of Tianjin under grant (No. 20JCYBJC01410).

\end{document}